\documentclass[a4paper,11pt]{article}

\usepackage{latexsym}
\usepackage{amsmath}
\usepackage{amsthm}

\usepackage{graphicx}
\usepackage{txfonts}
\usepackage{bm}
\usepackage{color}
\usepackage{epic,eepic}

\usepackage{geometry}
\numberwithin{equation}{section}

\newtheorem{theorem}{Theorem}[section]
\newtheorem{corollary}[theorem]{Corollary}
\newtheorem{lemma}[theorem]{Lemma}
\newtheorem{claim}[theorem]{Claim}
\newtheorem{remark}{Remark}[section]
\newtheorem{algorithm}[theorem]{Algorithm}

\newcommand{\OMIT}[1]{{\bf [OMIT:} #1 \ {\bf --- end OMIT] }}  
   \renewcommand{\OMIT}[1]{}            

\newcommand{\ZZ}{{\bf Z}}

\newcommand{\finbox}{\hspace*{\fill}$\rule{0.17cm}{0.17cm}$}
\newcommand{\finboxHere}{\ $\rule{0.17cm}{0.17cm}$}

\newcommand{\smallbox}{\scalebox{0.6}{\mbox{$\square$}}}
\newcommand{\sq}{\sp{\smallbox}}
\newcommand{\odotZ}{\overset{....}} 

\newcommand{\llceil}{\bigg\lceil} 
\newcommand{\rrceil}{\bigg\rceil}

\newcommand{\Proof}{\noindent {\bf Proof.  }}

\usepackage{lineno}
\newcommand*\patchAmsMathEnvironmentForLineno[1]{
  \expandafter\let\csname old#1\expandafter\endcsname\csname #1\endcsname
  \expandafter\let\csname oldend#1\expandafter\endcsname\csname end#1\endcsname
  \renewenvironment{#1}
     {\linenomath\csname old#1\endcsname}
     {\csname oldend#1\endcsname\endlinenomath}}
\newcommand*\patchBothAmsMathEnvironmentsForLineno[1]{
  \patchAmsMathEnvironmentForLineno{#1}
  \patchAmsMathEnvironmentForLineno{#1*}}
\AtBeginDocument{
\patchBothAmsMathEnvironmentsForLineno{equation}
\patchBothAmsMathEnvironmentsForLineno{align}
\patchBothAmsMathEnvironmentsForLineno{flalign}
\patchBothAmsMathEnvironmentsForLineno{alignat}
\patchBothAmsMathEnvironmentsForLineno{gather}
\patchBothAmsMathEnvironmentsForLineno{multline}
}

\begin{document}

\title{
Decreasing Minimization on M-convex Sets:  \\
Algorithms and Applications
}

\author{Andr\'as Frank\thanks{MTA-ELTE Egerv\'ary Research Group,
Department of Operations Research, E\"otv\"os University, P\'azm\'any
P. s. 1/c, Budapest, Hungary, H-1117. 
e-mail: {\tt frank@cs.elte.hu}. 
ORCID: 0000-0001-6161-4848.
The research was partially supported by the
National Research, Development and Innovation Fund of Hungary
(FK\_18) -- No. NKFI-128673.
}
\ \ and \ 
{Kazuo Murota\thanks{Department of Economics and Business Administration,
Tokyo Metropolitan University, Tokyo 192-0397, Japan, 
e-mail: {\tt murota@tmu.ac.jp}. 
Currently at
The Institute of Statistical Mathematics,
Tokyo 190-8562, Japan.
ORCID: 0000-0003-1518-9152.
The research was supported by 
CREST, JST, Grant Number JPMJCR14D2, Japan, and 
JSPS KAKENHI Grant Numbers JP26280004, JP20K11697. 
}}}

\date{July 2020 / July 2021}

\maketitle

\begin{abstract} 
This paper is concerned with algorithms and applications 
of decreasing minimization on an M-convex set, 
which is the set of integral elements of an integral base-polyhedron.  
Based on a recent characterization of decreasingly minimal (dec-min) elements,
we develop a strongly polynomial algorithm for computing a dec-min element
of an M-convex set.  
The matroidal feature of the set of dec-min elements makes it 
possible to compute a minimum cost dec-min element, as well.  
Our second goal is to exhibit various applications in matroid and network
optimization, resource allocation, and (hyper)graph orientation.  
We extend earlier results on semi-matchings 
to a large degree
by developing a structural description 
of dec-min in-degree bounded orientations of a graph.  
This characterization gives rise to a strongly polynomial algorithm 
for finding a minimum edge-cost dec-min orientation.
\end{abstract}


M-convex set, 
Polynomial algorithm.



{\bf Running head}:
Decreasing Minimization: Algorithms and Applications


\newpage
\tableofcontents

\newpage


\section{Introduction}
\label{SCintroB}

This paper is concerned with algorithms and applications of
decreasing minimization on an M-convex set, which is the set of
integral elements of an integral base-polyhedron.  
An element of a set of vectors, in general, 
is called {\bf decreasingly minimal} (dec-min)
if its largest component is as small as possible, within this, 
its second largest component is as small as possible, and so on.
Decreasing minimization means the problem of finding 
a dec-min element of a given set of vectors 
(or even a cheapest dec-min element with respect to a given linear cost-function).  
When the given set of vectors consists of integral vectors, 
this problem is also referred to as discrete decreasing minimization.  
In the literature, typically the term lexicographic 
optimization is used, 
but we prefer ``decreasing minimization''
because we also consider its natural counterpart 
``increasing maximization,'' and the use of these two symmetric terms
seems more appropriate to distinguish the two related notions.  
An element of a set of vectors is called {\bf increasingly maximal}
(inc-max) if its smallest component is as large as possible, within
this, its second smallest component is as large as possible, and so on.

In the companion paper \cite{FM19partA}, 
the present authors have investigated the structural aspects of 
the discrete decreasing minimization on an M-convex set.
Among others, the dec-min elements are characterized as those admitting no local improvement,
where the precise meaning of a local improvement is
defined formally in term of ``1-tightening step'' in Section~\ref{SCalgo2}. 
It is also shown that an element of an M-convex set is
decreasingly minimal precisely if it is a minimizer of 
the sum of the squared components.

As dual objects to dec-min elements, the notions of canonical chain, 
canonical partition of the ground-set, and 
essential-value sequence were defined, 
and the structure of the set of all dec-min elements 
is described in terms of these dual objects.  
(An equivalent definition of these notions
is given in Section \ref{SCcompcanchain}.)
We emphasize that the role of these dual objects 
is not merely to help us
fully understand the problem from its dual side.  
Beyond this, the dual characterization reveals the fundamental feature of the primal problem 
that the set of dec-min elements itself forms an M-convex set, and, in fact, a rather special
one arising from a matroid by translation.  
In addition, these dual objects are inherent in computing a dec-min element 
in strongly polynomial time and indispensable 
for efficient computation of a minimum weight dec-min element, as well.

The first goal of this paper is to develop,
on the basis of the above-mentioned structural characterizations, 
a strongly polynomial algorithm 
for computing
a dec-min element as well as  the canonical chain of a given M-convex set.  
The second goal is to exhibit several applications.  
For example, we prove a conjecture of Borradaile et al.~\cite{BIMOWZ} 
on dec-min strongly connected orientations of undirected graphs.  
Our general approach makes it possible to solve algorithmically 
even the minimum edge-cost dec-min orientation problem
when upper and lower bounds are imposed on the in-degrees and the
orientation is expected to be $k$-edge-connected (or even $(k,\ell)$-edge-connected).  
These orientation results form the basis of a major generalization of 
the so-called semi-matching problem initiated by 
Harvey et al.~\cite{HLLT}, 
which had been motivated by a resource allocation problem
in computer science.
Our approach is the first one that provides a strongly polynomial algorithm 
for the capacitated case, as well.

An algorithmic solution to a discrete counterpart of Megiddo's 
lexicographic flow problem  
\cite{Megiddo74,Megiddo77} is also developed.  
Yet another application of the structural results of \cite{FM19partA} gives
rise to an extension of a result of Levin and Onn \cite{Levin-Onn} on finding $k$
bases of a matroid on a ground-set $S$ with $n$ elements such that the
degree-vector of the hypergraph formed by these $k$ bases is
decreasingly minimal.  Our approach generalizes this problem to the
case when one has $k$ distinct matroids on $S$.

The paper is organized as follows.
Algorithms for computing
a dec-min element and the canonical chain
are given in Section \ref{SCalgo2}. 
In Section~\ref{alk}, various kinds of applications are shown, 
including those to matroids, network flows, arborescences, and connectivity augmentations. 
Sections \ref{SCori1}, \ref{SCori2}, and \ref{SCori3}
are devoted to detailed account of applications to graph orientation problems.

\subsection{Notation and terminology}
\label{SCnotation}

We continue to use notation and terminology introduced in \cite{FM19partA}, 
while some additional ones are given here.  
Let $S$ be a finite ground-set $S$ with $n$ elements.
Two subsets $X$ and $Y$ of $S$ are {\bf intersecting} 
if $X\cap Y\not =\emptyset$, 
{\bf properly intersecting} if none of $X\cap Y$, $X-Y$, and $Y-X$ is empty, 
and {\bf crossing} if none of $X-Y$, $Y-X$, $X\cap Y$, and $S-(X\cup Y)$ is empty.
For distinct elements $s, t$ of $S$, a subset of $S$ containing $t$
but not $s$ is called a {\bf $t\overline{s}$-set}.

For a vector $x\in {\bf R} \sp{S}$ or a function $x:S\rightarrow {\bf R}$, 
we define the set-function 
$\widetilde x:2\sp{S}\rightarrow {\bf R}$ 
by 
$\widetilde x(Z):=\sum [x(s):s\in Z]$ $(Z\subseteq S)$. 
This is a modular function in the sense that 
$\widetilde x(X)+\widetilde x(Y) = \widetilde x(X\cap Y)+\widetilde x(X\cup Y)$ 
holds for every $X,Y\subseteq S$.

For any integral polyhedron 
$R\subseteq {\bf R} \sp{S}$, 
we use the notation $\odotZ{R}$ to denote the set of integral elements of $R$,
that is, 
\begin{equation} \label{odotBdef} 
\odotZ{R} := R\cap {\bf Z}\sp{S}, 
\end{equation}
where $\odotZ{R}$ may be pronounced  ``dotted $R$.''  
The notation is intended to refer intuitively to the set of lattice points of $R$.

For a set-function $h$, we allow it to have value $+\infty $ or $-\infty $. 
Unless otherwise stated, $h(\emptyset )=0$ is assumed throughout.  
When $h(S)$ is finite, the {\bf complementary function} $\overline{h}$ 
is defined by 
$\overline{h}(X):  = h(S) - h(S-X)$.  
Observe that the complementary function of $\overline{h}$ is $h$ itself.

Let $b$ be a set-function with $b(\emptyset )=0$,
for which $b(X)=+\infty $ is allowed but $b(X)=-\infty $ is not.  
The submodular inequality for subsets $X,Y\subseteq S$ is defined by 
\begin{equation} \label{subfndef} 
b(X) + b(Y) \geq b(X\cap Y) + b(X\cup Y).
\end{equation} 
We say that $b$ is (fully) {\bf submodular} if this inequality holds 
for every pair of subsets $X, Y\subseteq S$ with finite $b$-values.  
When the submodular inequality is required only
for intersecting (crossing) pairs of subsets, 
we say that $b$ is {\bf intersecting} ({\bf crossing}) {\bf submodular}.  
A set-function $p$ is called (fully, intersecting, crossing) {\bf supermodular} 
if $-p$ is (fully, intersecting, crossing) submodular.  
It follows from the definitions that 
if $b$ is fully (or crossing) submodular with finite $b(S)$, 
then its complementary function $\overline{b}$ is 
a fully (or crossing) supermodular function, 
while $\overline{b}$ is not necessarily intersecting supermodular 
when $b$ is intersecting submodular.

For a (fully) submodular integer-valued set-function $b$ on $S$  
with $b(\emptyset )=0$ and finite $b(S)$, 
the {\bf base-polyhedron} $B$ is defined by 
\begin{equation} \label{Bsubfn} 
 B=B(b) = \{x\in {\bf R}\sp{S}:  \widetilde x(S)=b(S), 
     \ \widetilde x(Z)\leq b(Z) \ \mbox{ for every } \ Z\subset S\}, 
\end{equation}
which is a (possibly unbounded) integral polyhedron in ${\bf R}\sp{S}$.  
Section~14 of book \cite{Frank-book} provides an overview
of basic properties of base-polyhedra;
see also the book of Schrijver \cite{Schrijverbook}.  
For example, it is a basic property that $B=B(b)$ is a non-empty integral polyhedron, 
and $B$ uniquely determines its defining (fully) submodular function $b$, namely, 
\[
b(Z)= \max \{\widetilde x(Z):  x\in B\} 
\quad ( = \max \{\widetilde x(Z):  x\in \odotZ{B}\}).
\]
By convention, the empty set is also considered a base-polyhedron.
A (fully) supermodular integer-valued set-function $p$ with
$p(\emptyset )=0$ and $p(S)$ finite also defines an integral base-polyhedron by 
\begin{equation} \label{Bsupfn} B=B'(p)=\{x \in
{\bf R}\sp{S}:  \widetilde x(S)=p(S), 
 \ \widetilde x(Z)\geq p(Z) \  \mbox{ for every } \ Z\subset S\} ,
\end{equation} 
since the complementary function $b:= \overline{p}$ of $p$ 
is fully submodular and $B'(p) = B(b)$. 
In discrete convex analysis \cite{Murota98a,Murota03}, the set $\odotZ{B}$ 
of integral elements of an integral base-polyhedron $B$ is called an {\bf M-convex set}.

We say that a submodular function $b$ and a supermodular function $p$ meet 
the {\bf cross-inequality} for $X,Y\subseteq S$ if 
\begin{equation} \label{crossineqdef} 
b(X)-p(Y) \geq b(X-Y) - p(Y-X).  
\end{equation}
A pair $(p,b)$ of set-functions is called {\bf fully paramodular} 
(or we say that $(p,b)$ is a {\bf strong pair}) 
if $b$ is fully submodular, $p$ is fully supermodular, and 
the cross-inequality holds for every $X,Y\subseteq S$.  
For a strong pair $(p,b)$, the polyhedron $Q$ defined by 
\begin{equation} \label{QpbGpoly} 
Q=Q(p,b) := \{x \in {\bf R}\sp{S}:  
  p(Z) \leq \widetilde x(Z) \leq b(Z) \ \mbox{ for every } \ Z\subseteq S \} 
\end{equation} 
is called a {\bf generalized polymatroid} 
({\bf g-polymatroid}, for short), 
while $(p,b)$ is called the {\bf border pair} of $Q$.  
It is a basic fact \cite{Frank-book} that $Q$ is never empty.  
By convention, the empty set is also considered a g-polymatroid.  
Furthermore, a non-empty g-polymatroid $Q$ uniquely determines 
its fully paramodular border pair, and $Q$ is integral when
$p$ and $b$ are integer-valued.  
Base-polyhedra are special g-polymatroids, namely, those for which $p(S)=b(S)$, 
and every g-polymatroid arises from a base-polyhedron 
by projection along a single coordinate axis.

To prove theorems on base-polyhedra, 
it is much easier to work with base-polyhedra defined 
by fully sub- or supermodular functions.  
For applications, however, it is fundamentally important 
that weaker set-functions
(e.g. intersecting, crossing or even weaker submodular functions) 
may also define base-polyhedra (or M-convex sets), as well as g-polymatroids.  
(We shall use this fact frequently in Sections \ref{SCori2} and \ref{SCori3}.)  
For example, if $p$ is an integer-valued intersecting or crossing supermodular function, 
then $B'(p)$ can be proved to be an integral base-polyhedron 
(see, e.g. Theorem~15.3.4 in book \cite{Frank-book}), 
which may, however, be empty.  
(This result, for example, underlies the fact
that the in-degree vectors of $k$-edge-connected orientations of a
$2k$-edge-connected graph form an M-convex set.)

When $p$ is intersecting supermodular, $b$ is
intersecting submodular, and the cross-inequality holds for properly
intersecting pairs of subsets $X,Y$, we say that 
$(p,b)$ is a {\bf weak pair} or {\bf intersecting paramodular}.  
For such a border pair, $Q=Q(p,b)$ is known to be a (possibly empty) g-polymatroid.  
It should be noted, that the (unique) fully paramodular border pair defining $Q$
can be concretely expressed with the help of the weak pair $(p,b)$,
but this formula is quite complicated 
(see, Corollary 15.3.4 in book \cite{Frank-book}).

It is a non-trivial task to characterize the situation when $B(b)$ is non-empty 
for a crossing submodular function $b$
but Fujishige \cite{Fujishige84e} developed 
an elegant necessary and sufficient condition.  
A similar theorem can be formulated for
g-polymatroids defined by a weak pair $(p,b)$ 
(see, Theorem 15.3.13 in book \cite{Frank-book}).  
It should, however, be emphasized that applications 
often need base-polyhedra or g-polymatroids 
defined by even weaker sub- or supermodular functions.  
For example, the in-degree vectors of 
$k$-edge-connected digraphs
obtained by adding a given number of arcs
to an input digraph $H=(V,A)$ form an M-convex set.  
(See, Theorem~17.2.9 and Section~15 of \cite{Frank-book} for an overview.  
A recent paper of B{\'e}rczi and Frank \cite{FrankJ67} 
includes even more intricate constructions of M-convex sets 
appearing in graph connectivity augmentation problems.)

We assume that graphs or digraphs have no loops but parallel edges are allowed.  
Sometimes we refer to an edge of a digraph as an {\bf arc}.
For a digraph $D=(V,A)$, the {\bf in-degree} of a node $v$ is 
the number of arcs of $D$ with head $v$.  
The in-degree $\varrho_{D}(Z)=\varrho(Z)$ of a subset $Z\subseteq V$ 
denotes the number of edges ($=$ arcs) entering $Z$, 
where an arc $uv$ is said to enter $Z$ if its head $v$ is in $Z$ 
while its tail $u$ is in $V-Z$.
The {\bf out-degree} $\delta_{D}(Z)=\delta(Z)$ is the number of arcs leaving $Z$,
 that is $\delta(Z)=\varrho(V-Z)$.  
The number of edges of a directed or undirected graph $H$ 
induced by $Z\subseteq V$ is
denoted by $i(Z)=i_H(Z)$.  
In an undirected graph $G=(V,E)$, the {\bf degree} $d(Z)=d_{G}(Z)$ 
of a subset $Z\subseteq V$ denotes the number of edges  connecting $Z$ and $V-Z$,
 while $e(Z)=e_{G}(Z)$ denotes the number of edges with one or two end-nodes in $Z$.  
Clearly, $e(Z)=d(Z)+i(Z)$.

\section{Algorithms} 
\label{SCalgo2}

In this section, we consider algorithmic aspects of
decreasing minimization over an  M-convex set.
In particular, we show how to compute efficiently 
a decreasingly minimal element along with 
its canonical chain and partition.

Let $B$ be a non-empty integral base-polyhedron and let $p$ denote the
unique fully supermodular function for which $B=B'(p)$.  
Let $m$ be an element of the M-convex set $\odotZ{B}$.  
We need some definitions introduced in \cite{FM19partA}.  
A {\bf 1-tightening step} replaces $m$ by $m':=m+\chi_{s}-\chi_{t}$,
where $s$ and $t$ are elements of $S$
for which $m(t)\geq m(s)+2$ and $m'$ belongs to $\odotZ{B}$.  
A subset $X\subseteq S$ is called {\bf $m$-tight}
(with respect to $p$)
if $\widetilde m(X)=p(X)$.  
A subset $X\subseteq S$ is called an {\bf $m$-top} set
if $m(u)\geq m(v)$ holds whenever $u\in X$ and $v\in S-X$.
We call an integral vector $x\in {\bf Z}\sp{S}$ 
{\bf near-uniform} on a subset $S'$ of $S$ if
its largest and smallest components on $S'$ differ by at most 1, that is, if
$x(s)\in \{\ell, \ell +1\}$ for some integer $\ell$ for every $s\in S'$.

First we recall fundamental characterizations of a dec-min element of an M-convex set.

\begin{theorem}[{\cite[Theorem 3.3]{FM19partA}}]  
\label{THequidotone}
For an element $m$ of an M-convex set $\odotZ{B} = \odotZ{B'}(p)$,
the following four conditions are pairwise equivalent.
\smallskip

\noindent
{\rm (A)} \ There is no 1-tightening step for $m$.

\smallskip

\noindent 
{\rm (B)} \ There is a chain 
$(\emptyset \subset) \ C_{1}\subset C_{2}\subset \cdots \subset C_\ell \ ( = S)$ 
such that 
each $C_{i}$ is an $m$-top and $m$-tight set (with respect to $p$) 
and $m$ is near-uniform on each 
$S_{i}:= C_{i}-C_{i-1}$ $(i=1,2,\dots ,\ell)$,
where $C_{0}:=\emptyset$.

\smallskip

\noindent 
{\rm (C1)} \ $m$ is decreasingly minimal in $\odotZ{B}$.
\smallskip

\noindent 
{\rm (C2)} \ $m$ is increasingly maximal in $\odotZ{B}$.  
\finbox
\end{theorem}

We mentioned already in Introduction that it is 
convenient to prove results for a base-polyhedron 
assuming  that it is given by a fully sub- or supermodular set-function.
In applications, however,
a base-polyhedron is typically given 
by an intersecting or crossing (or even weaker) function.
Therefore, in describing and analysing algorithms, 
we consider these weaker functions, as well.

\begin{remark} \rm \label{RMsubmmin}
One of the most fundamental algorithms of discrete optimization is 
for minimizing a submodular function, that is, 
for finding a subset $Z$ of $S$ for which 
$b(Z)=\min \{b(X):X\subseteq S\}$.  
There are strongly polynomial algorithms for this problem;
Schrijver \cite{Schrijver2000} and Iwata et al.~\cite{IFF01} are the first.  
We shall refer to such an algorithm as a {\bf submod-minimizer} subroutine.  
The complexity of a significantly more efficient algorithm, 
due to Orlin \cite{Orlin09}, is $O(n\sp{6})$ (where $n= |S|$)
and this algorithm calls $O(n\sp{5})$ times a routine which evaluates 
the submodular function in question.  
(An evaluation routine outputs the value $b(X)$ for any input subset $X\subseteq S$). 
 A recent algorithm of Jiang \cite{Jiang21} needs only $O(n\sp{3})$ calls
of an evaluation routine.
Naturally, submodular function minimization and supermodular function
maximization are equivalent.
\finbox
\end{remark}

\subsection{The basic algorithm for computing a dec-min element}
\label{SCbasicalg}

Our first goal is to describe a natural 
approach---the {\bf basic algorithm}---for finding a decreasingly minimal element 
of an M-convex set $\odotZ{B}$.  
The basic algorithm is always finite, and if $B\subseteq {\bf R}_{+}\sp{S}$, 
it is pseudopolynomial in the sense that it is polynomial in $n + |p(S)|$. 
This means that the algorithm is polynomial in $n$ when 
$B\subseteq {\bf R}_{+}\sp{S}$ and $|p(S)|$ can be bounded by a polynomial of $n$. 
This is the case, for example, 
in an application to
strongly connected decreasingly minimal
(=egalitarian) orientations.  
In the general case, where typical
applications arise by defining $p$ with a `large' capacity function, 
a (more sophisticated) 
strongly polynomial-time algorithm will be
described in Section \ref{strong.pol}.
In order to find a dec-min element of an M-convex set $\odotZ{B}$, 
we assume that a subroutine is available to 
\begin{equation} 
\hbox{compute an integral element of $B$.}\ 
\label{(member)} 
\end{equation}
We use this subroutine only once to obtain an initial point in the algorithm.

Suppose now that an integral member $m$ of $B$ is available.
The algorithm needs a subroutine to 
\begin{equation} 
\hbox{ decide for  $m\in \odotZ{B}$ and for $s,t\in S$ if $m':=m+\chi_{s}-\chi_{t}$ 
belongs to $B$.}\ 
\label{(st.step)} 
\end{equation}
By $n\sp 2$ applications of subroutine \eqref{(st.step)}, 
one can decide, for a given $m\in \odotZ{B}$, 
whether there exists a 1-tightening step or not.
Observe that the subroutine \eqref{(st.step)} is certainly available if we can
\begin{equation} 
\hbox{ decide for any $m'\in {\bf Z}\sp{S}$ whether or not $m'$ belongs to $B$,}\
\label{(st.step2)} 
\end{equation}
though applying this more general subroutine is clearly
slower than a direct algorithm to realize \eqref{(st.step)}.
In this paper we mainly work with Subroutines \eqref{(member)} and \eqref{(st.step)},
which are representation-free in the sense that they do not refer directly
to the supermodular (or submodular) functions describing $B$.
The realization of these subroutines depends on how 
the base-polyhedron $B$ (or an M-convex set $\odotZ{B}$) is given.
Recall that a base-polyhedron can be described by a fully,
intersecting, crossing, or even weaker functions.
For example, when $B$ is described by a fully supermodular function $p$,
\eqref{(member)} can be carried out with $n$ evaluations of $p$
and \eqref{(st.step)} can be done with a single application of a submod-minimizer.
See Remark~\ref{RMinsercross} for more details.

\medskip

As long as possible, apply the 1-tightening step.
By Theorem~\ref{THequidotone},
when no more 1-tightening step is available,
the current $m$ is a decreasingly minimal member of $\odotZ{B}$ 
and the algorithm terminates.  
In this connection we may recall  
the following characterization of a dec-min element of an M-convex set.

\begin{theorem}[{\cite[Corollary 6.4]{FM19partA}}]  
\label{THdecmin=squaresum}
An element $m$ of an M-convex set $\odotZ{B}$ is
a dec-min element of $\odotZ{B}$
if and only if it is a minimizer of
the square-sum 
$W(z) :=\sum [z(s)\sp{2}:s\in S]$ 
of $z$ over the elements $z$ of $\odotZ{B}$.
\finbox
\end{theorem}

Now observe that a single 1-tightening step strictly decreases the
square-sum of the components.
This implies that
the number of 1-tightening steps is
bounded by $W(m)$ when $m$ is the initial member of $\odotZ{B}$. 
In particular, the algorithm is finite.

This algorithm, however, may be quite inefficient as is demonstrated by
the simple example where $\vert S\vert =2$ and 
$B=\{(x_1,x_2):  \ x_1+x_2 = 0\}$.  
Here if the initial member of $\odotZ{B}$ is, for example, 
$m=(10\sp{6}, - 10\sp{6})$, then the algorithm needs 
$10\sp{6}$ 1-tightening steps.  
However, if $B$ is in the non-negative orthant
(which is often the case in applications), then the following
reasonable bound can be given for the complexity.
Since the square-sum of an arbitrary integral vector $z\geq 0$
with $\widetilde z(S)=p(S)$ is at most
$p(S)\sp{2}$ and $\widetilde z(S)=p(S)$
holds for all members $z$ of 
$\odotZ{B}$,
we conclude that the number of 1-tightening steps 
is at most $p(S)\sp{2}$.  
Therefore, if $B \subseteq {\bf R}_{+}\sp{S}$ and
$\vert p(S)\vert $ is bounded by a polynomial of $n$, 
then the basic algorithm to compute a dec-min
element of $\odotZ{B}$ is strongly polynomial.

Although the basic algorithm is efficient when $B\subseteq {\bf R}_{+}\sp{S}$ 
and $|p(S)|$ is `small'  
(that is, $\vert p(S)\vert $ is bounded by a power of $n$), 
it is not strongly polynomial when $|p(S)|$ is `large'.  
We postpone, till Section \ref{strong.pol},
the description of a strongly polynomial algorithm for computing a
dec-min element of an M-convex set $\odotZ{B}$ defined by a general $p$.  
In the next section we show how the canonical chain as well as
the essential value-sequence can be computed, once a dec-min element $m$ is available.  
It is emphasized that these dual objects are
indispensable and must be computed when we are interested 
in identifying the set of all dec-min elements of $\odotZ{B}$ or 
in finding a minimum weight dec-min element (cf., \cite[Section~5.3]{FM19partA}).

\begin{remark} \rm \label{RMinsercross}
In the algorithm above, we did not rely
explicitly on the set-function defining the base-polyhedron $B$ in
question, apart from the single value $p(S)$.  
The only assumption was
that $B$ is non-empty and the subroutines 
\eqref{(member)} and \eqref{(st.step)} are available.  
When the base-polyhedron $B=B'(p)$ is given by a fully supermodular function $p$, 
then the subroutine \eqref{(member)} 
 can be realized by a version of
the polymatroid greedy algorithm of Edmonds \cite{Edmonds70} that
needs $n$ evaluations of $p$ (but not a submod-minimizer subroutine).  
If $p$ is intersecting supermodular, then an algorithm of 
Frank and Tardos \cite{FrankJ17} requires $n$ calls of a submod-minimizer.  
The same paper includes an algorithm for the case when $p$ is crossing
supermodular, and this algorithm requires $n\sp{2}$ calls of a submod-minimizer.
Concerning the other subroutine
\eqref{(st.step)}, we note that  $m'=m+\chi_{s}-\chi_{t}$ is in $B=B'(p)$ precisely 
if there is no $m$-tight 
$t\overline{s}$-set
(with respect to $p$),
and this is true even when $B$ is defined by a crossing supermodular function $p$.  
Therefore, 
the subroutine \eqref{(st.step)}
can be realized with a single call of a submod-minimizer 
if $p$ is fully supermodular, 
$n$ calls if  $p$ is intersecting supermodular, and
$n\sp{2}$ calls if  $p$ is crossing supermodular.
\finbox
\end{remark}

\begin{remark} \rm \label{RMsqsumminalg}
Theorem \ref{THdecmin=squaresum} shows that a dec-min element of an M-convex set
is characterized as a minimizer of the square-sum $W(z)$.
This naturally suggests the approach of minimizing $W(z)$ to find a dec-min element.
Minimizing the square-sum $W(z)$ over an M-convex set
is a special case of minimizing an M-convex function,
for which a local improvement algorithm works
(in finite steps) \cite{Murota03}.
The basic algorithm described above corresponds to the special case of 
this local improvement algorithm for M-convex function minimization.
It is also noted that
Fujishige \cite{Fujishige80}
(see also \cite[Section 8.2]{Fujishigebook})
solved the continuous case of the square-sum minimization in polynomial time,
which problem is equivalent to finding 
the unique minimum norm point $m\sp{*}$ of a base-polyhedron $B$.
In \cite{FM19partII}
 we proved a theorem formalizing the  intuitive feeling that all  the
 dec-min elements of $\odotZ{B}$ are in the neighbourhood of 
$m\sp{*}$, and this
 result makes it possible to develop an alternative algorithm to compute a
 dec-min element of  $\odotZ{B}$. 
\finbox
\end{remark}

\subsection{Computing the essential value-sequence and the canonical chain} 
\label{SCcompcanchain}

In this section, we describe an algorithm that assigns a chain of
subsets of $S$, a partition of $S$, and a strictly decreasing sequence
of integers to a given dec-min element $m$ of an M-convex set.

Let $B=B'(p)$ be again a (non-empty) integral base-polyhedron whose unique 
(fully) supermodular bounding function is $p$.  
In the algorithm,  
we assume that we can  compute the smallest $m$-tight set 
$T_{m}(u) = T_{m}(u;p)$ 
containing a given element $u\in S$.
Since we have
$T_{m}(u)=\{s:  m+\chi_{s} -\chi_{u}\in B\}$,
$T_{m}(u)$ is indeed computable by at most $n$ applications 
of Subroutine \eqref{(st.step)}.

\begin{algorithm}  \rm  \label{algo1}
Given a dec-min element $m$ of $\odotZ{B}$, 
the following procedure computes 
a chain
${\cal C}\sp{*}=\{C_{1},C_{2},\dots ,C_{q}\}$
with $C_{1} \subset C_{2} \subset \cdots \subset C_{q} \ (=S)$ 
and a partition 
${\cal P}\sp{*}=\{S_{1},S_{2}, \dots , S_{q}\}$ of $S$
along with a sequence
$\beta_{1}>\beta_{2}>\cdots >\beta_{q}$ of integers.

\begin{enumerate}
\item 
 Let $\beta_{1}$ denote the largest value of $m$.
Let $C_{1} := \bigcup \{  T_{m}(u):  m(u) =\beta_{1} \}$,
$S_{1} := C_{1}$, and $i := 2$.

\item 
In the general step $i \geq 2$, 
the pairwise disjoint non-empty sets $S_{1},S_{2},\dots ,S_{i-1}$ 
and a chain $C_{1} \subset C_{2} \subset \dots \subset C_{i-1}$
have already been computed along with 
the values 
$\beta_{1}>\beta_{2}> \dots > \beta_{i-1}$.
If $C_{i-1} = S$, set $q:= i-1$ and stop.
Otherwise, let 
\begin{align*}
\beta_{i} & := \max\{m(s):  s\in S-C_{i-1} \},
\\
C_{i}  &  := \bigcup \{ T_{m}(u):  m(u)\geq \beta_{i}\}, 
\\
S_{i}   & :=C_{i}-C_{i-1},
\end{align*}
and go to the next step for $i:=i+1$.
\finbox
\end{enumerate}
\end{algorithm}

It was proved in \cite{FM19partA} (Corollary~5.4) 
that these sequences do not depend on the choice of $m$.  
Therefore the chain ${\cal C}\sp{*}$ is called the 
{\bf canonical chain} belonging to $\odotZ{B}$, 
the partition ${\cal P}\sp{*}$ is the {\bf canonical partition} of $\odotZ{B}$, 
while the sequence $\{\beta_{1},\beta_{2}, \dots ,\beta_{q} \}$ 
is called the {\bf essential value-sequence} of $\odotZ{B}$.  
We emphasize that Algorithm \ref{algo1}
is strongly polynomial for arbitrary $p$
(independently of the magnitude of $\vert p(S)\vert $), provided that
a dec-min element $m$ of $\odotZ{B}$ is already available 
as well as the subroutine \eqref{(st.step)}.

It is in order here to emphasize
the significance of this  algorithm for computing these dual objects.
By Theorem \ref{THchardecmin} below,
Algorithm \ref{algo1} enables us to computationally capture 
the set of all dec-min elements.
Concisely, the matroid associated with dec-min elements,
as in Theorem \ref{THmatroid-eltolt} below,
can be identified by this algorithm.
Here a {\bf matroidal M-convex set} means \cite{FM19partA}
the translation of the incidence vectors 
of bases of a matroid by an integral vector.

\begin{theorem}[{\cite[Corollary 5.2]{FM19partA}}] \label{THchardecmin}
An element $m$ of an M-convex set $\odotZ{B}$ is decreasingly minimal 
if and only if each 
$C_{i} \in {\cal C}\sp{*}$ 
is $m$-tight (with respect to $p$) 
and $\beta_{i} -1 \leq m(s) \leq \beta_{i}$ holds for each 
$s\in S_{i}$ \ $(i=1,\dots ,q)$.  
\finbox 
\end{theorem}

\begin{theorem}[{\cite[Theorem 5.7]{FM19partA}}]  \label{THmatroid-eltolt} 
The set of dec-min elements of an M-convex set $\odotZ{B}$ is a matroidal M-convex set.
\finbox 
\end{theorem} 

We shall we use Theorem~\ref{THchardecmin} in Sections \ref{SCoricano} and \ref{SCdecminfg},
and Theorem~\ref{THmatroid-eltolt} in Sections \ref{inout} and \ref{SCoricheapest}.

\paragraph{Adaptation to the intersection with a box} \ 
Algorithm \ref{algo1} can be adapted to the case when we have specific upper
and lower bounds on the members of $\odotZ{B}= \odotZ{B'}(p)$.  
Let
$f:S\rightarrow {\bf Z} \cup \{-\infty \}$ and 
$g:S\rightarrow {\bf Z} \cup \{+\infty \}$ 
be bounding functions with $f\leq g$ and let 
$T(f,g):=\{x\in {\bf R}\sp{S}:  f\leq x\leq g\}$
denote the box defined by $f$ and $g$.  
It is a basic fact that the intersection of an integral base-polyhedron 
with an integral box is another integral base-polyhedron
(see, Theorem 14.3.9 in book \cite{Frank-book} for this statement 
in a more general context).
Therefore the intersection $B\sq :=B\cap T(f,g)$ is also a
(possibly empty) integral base-polyhedron.  
Assume that ${B\sq}$ is non-empty.

Let $m$ be an element of $\odotZ{B\sq}$
($= B\sq \cap \ZZ\sp{S}$).
Let $T_{m}(u)$ denote the smallest $m$-tight set containing $u$ with respect to $p$, 
and let $T\sq_{m}(u)$ be the smallest $m$-tight set containing $u$ with respect to $p\sq$.

\begin{claim} \label{CLTm(u)}
\[
T\sq_{m}(u) = 
\begin{cases} \{u\} & 
\ \ \hbox{\rm if}\ \ \ m(u)=f(u),
\\ 
T_{m}(u) - \{v:  m(v) = g(v)\} & \ \ \hbox{\rm if}\ \ \ m(u)>f(u).
\end{cases}
\]
\end{claim}  

\Proof 
We have $T\sq_{m}(u)= \{ s: m-\chi_{u} + \chi_{s} \in B\sq \}$.
Since
$B\sq = B \cap T(f,g)$,
we have
$m-\chi_{u} + \chi_{s} \in B\sq$
if and only if
(i)
$m-\chi_{u} + \chi_{s} \in B$ and (ii) $m-\chi_{u} + \chi_{s} \in T(f,g)$ hold.
For $s \not = u$,  (i) holds 
if and only if
$s \in T_{m}(u)$,
and (ii) holds 
if and only if
$m(u) > f(u)$ and $m(s) < g(s)$.
Hence follows the claim.
\finbox \medskip

The claim implies that Algorithm \ref{algo1}
can be adapted easily to compute the
canonical chain and partition belonging to $\odotZ{B\sq}$
along with its essential value-sequence.

Our next goal is to describe a strongly polynomial algorithm to
compute a dec-min element of $\odotZ{B}$ in the general case when no
restriction is imposed on the magnitude of $\vert p(S)\vert $. 
To this end, we need an algorithm to maximize 
$\lceil {p(X) / \vert X\vert} \rceil$,
which is given in Section \ref{pbmaximizing}.
The strongly polynomial algorithm for computing a dec-min element
is described in Section \ref{strong.pol}.

\subsection{Maximizing $\lceil {p(X) / |X| }\rceil$ with the Newton--Dinkelbach algorithm} 
\label{pbmaximizing}

In this section we describe a variant
of the Newton--Dinkelbach (ND) algorithm to compute the maximum of
$\lceil {p(X) / |X|}\rceil $.
We assume that $p$ is an integer-valued set-function 
on a ground-set $S$ with $n\geq 1$ elements,
$p(\emptyset )=0$, and $p(S)$ is finite
($p(X)$ may be $-\infty$ for some $X$ but never $+\infty$).

An excellent overview by Radzik \cite{Radzik} 
presents fundamental properties of the ND-algorithm,
describing (among others) a strongly polynomial algorithm
for minimizing (or maximizing) the ratio of two modular set-functions.
For the problem of maximizing $\lceil {p(X) / |X|}\rceil$
(or ${p(X) / |X|}$), 
the ND-algorithm terminates in at most $n$ iterations,
which follows from the observation of Topkis \cite{Top78}
that the function 
$h(\mu) := \max \{ p(X) - \mu |X| : X  \subseteq S \}$
is a convex, piecewise-linear function with at most $n$ breakpoints
A recent paper by Goemans et al.~\cite {GGJ17}
establishes a quadratic 
$O(n\sp{2})$ bound on the number of iterations
of the ND-algorithm for maximizing $p(X) / a(X)$ over $X$ with $a(X)>0$,
where $a$ is an arbitrary modular set function.  
We present a variant of the ND-algorithm whose specific feature is that it works
throughout with integers $\lceil {p(X) / |X| }\rceil$. 
This has the advantage that the proof is simpler than the original one 
working with the fractions ${p(X) / |X| }$.

The algorithm works if a subroutine is available to 
\begin{equation} 
\hbox{ find a subset \, $X\subseteq S$ \, maximizing
\,  $p(X) - \mu |X|$ \, for any fixed integer $\mu$. } 
\label{(ND.routine)} 
\end{equation}
\noindent
This subroutine will actually be needed only for special values of $\mu $ when 
$\mu =\lceil p(X)/\ell\rceil$ 
(where $X\subseteq S$ and $1\leq \ell\leq n$). 
We do not have to assume that $p$ is supermodular and 
the only requirement for the ND-algorithm is that 
Subroutine \eqref{(ND.routine)} be available.  
Via a submod-minimizer this is certainly 
the case when $p$ happens to be
supermodular 
(cf., Remark \ref{RMsubmmin}).

In several applications, the requested general purpose
submod-minimizer can be superseded by a direct and more efficient
algorithm such as the one for network flows or for matroid partition.
The subroutine \eqref{(ND.routine)} is also available in the
more general case (needed in applications)
when $p$ is only crossing supermodular.  
Indeed, for a given ordered pair of elements $s,t\in S$, 
the restriction of $p$ on the family of $s\overline{t}$-sets
 is fully supermodular, 
and therefore we can apply a submod-minimizer 
to each of the $n(n-1)$ ordered pairs $(s,t)$ to get the requested maximum of $p(X) - \mu |X|$.

We call a value $\mu $ \ {\bf good} \ if $\mu |X| \geq p(X)$ 
for every $X\subseteq S$.  
A value that is not good is called {\bf bad}.  
Clearly, a sufficiently large $\mu $ is good.  
Our goal is to compute the minimum $\mu_{\rm min}$ of the good integers.  
This number is nothing but the maximum of  \ $\lceil {p(X) / |X| }\rceil $ \
over non-empty subsets of $S$.

Let $\mu_{0} := \lceil {p(S) / \vert S\vert } \rceil -1 $.
This (possibly negative) number is bad and the algorithm starts with $\mu _{0}$.  
Let 
\[
  X_{0} \in \arg\max \{ p(X)-\mu_{0}|X|  :  \ X\subseteq S \},
\]
that is, \ $X_{0}$ is a set maximizing the function $p(X)-\mu_{0}|X|$.  
Note that the badness of $\mu_{0}$ implies that 
$p(X_{0}) > \mu_{0}|X_{0}|$.

The procedure determines one by one a series of pairs $(\mu_{j},X_{j})$
for subscripts $j=1,2,\dots $ where each integer $\mu_{j}$ 
is a tentative candidate 
for $\mu $ while $X_{j}$ is a non-empty subset of $S$.
Suppose that the pair $(\mu_{j-1},X_{j-1})$
has already been determined for a subscript $j\geq 1$.  
Let $\mu_{j}$ be the smallest integer 
for which $\mu_{j} |X_{j-1}| \geq p(X_{j-1})$,
that is,
\[
 \mu_{j}:  =\llceil {p(X_{j-1}) \over |X_{j-1}|}\rrceil . 
\]

If $\mu_{j}$ is bad, that is, if there is a set $X\subseteq S$ with
$p(X) -\mu_{j}|X|  > 0$, then let 
\[
  X_{j} \in \arg\max \{ p(X)-\mu_{j}|X|  :  \ X\subseteq S \},
\] 
that is, \ $X_{j}$ is a set
maximizing the function \ $p(X)-\mu_{j}|X| $.  
(If there are more than one maximizing set, we can take any).  
Since $\mu_{j}$ is bad, we have
$X_{j}\not =\emptyset $ and $p(X_{j}) - \mu_{j}|X_{j}|>0$.

\begin{claim}  \label{betano} 
If $\mu_{j}$ is bad for some subscript $j\geq 0$, then $\mu_{j} < \mu_{j+1}$.  
\end{claim}

\Proof 
The badness of $\mu_{j}$ means that 
$p(X_{j})-\mu_{j}|X_{j}| > 0$, from which
\[
\mu_{j+1} = \llceil {p(X_{j}) \over |X_{j}|}\rrceil 
  = \llceil {p(X_{j})-\mu_{j}|X_{j}| \over |X_{j}| }\rrceil + \mu_{j} 
\ > \  \mu_{j}.  
\]
\vspace{-2.8\baselineskip} \\
\finbox
\vspace{1.5\baselineskip}

Since there is a good $\mu $ and the sequence $\mu_{j}$ is 
strictly monotone increasing by Claim \ref{betano}, 
there will be a first subscript $h\geq 1$ for which $\mu_{h}$ is good.  
The algorithm terminates by outputting this $\mu_{h}$ 
(and in this case $X_{h}$ is not computed).

\begin{theorem}   \label{Msteps} 
If $h$ is the first subscript during the run of the algorithm 
for which $\mu_{h}$ is good, then 
$\mu_{\rm min}=\mu_{h}$ 
(that is, $\mu_{h}$ is the requested smallest good $\mu $-value)
and $h\leq n$.
\end{theorem}

\Proof 
Since $\mu_{h}$ is good and $\mu_{h}$ is the smallest integer for
which $\mu_{h} |X_{h-1}| \geq p(X_{h-1})$, 
the set $X_{h-1}$ certifies that no good integer $\mu $ 
can exist which is smaller than $\mu_{h}$,
that is, $\mu_{\rm min}=\mu_{h}$.

\begin{claim}   \label{Xno} 
If $\mu_{j}$ is bad for some subscript $j\geq 1$, then $|X_{j-1}| > |X_{j}|$.  
\end{claim}

\Proof
As $\mu_{j}$ \ ($= \lceil {p(X_{j-1}) / |X_{j-1}|}\rceil$) \ is bad, 
we obtain that
\begin{align*}
p(X_{j})-\mu_{j}|X_{j}|  >0 
  & = p(X_{j-1}) - { p(X_{j-1}) \over |X_{j-1}|} |X_{j-1}| 
\\ &
\geq p(X_{j-1}) - \llceil {p(X_{j-1}) \over |X_{j-1}|}\rrceil |X_{j-1}|
= p(X_{j-1}) - \mu_{j}|X_{j-1}| ,
\end{align*}
from which we get
\begin{equation} \label{NDprfA}
 p(X_{j}) - \mu_{j}|X_{j}| > p(X_{j-1}) - \mu_{j}|X_{j-1}|.
\end{equation}
\noindent
Since $X_{j-1}$ maximizes $p(X) - \mu_{j-1}|X| $, we have
\begin{equation} \label{NDprfB}
p(X_{j-1}) - \mu_{j-1}|X_{j-1}| \geq p(X_{j}) - \mu_{j-1}|X_{j}|.  
\end{equation}
\noindent
By adding up 
\eqref{NDprfA} and \eqref{NDprfB},
we obtain
\[
 (\mu_{j} - \mu_{j-1})|X_{j-1}| > (\mu_{j} - \mu_{j-1})|X_{j}|.
\]
\noindent
As $\mu_{j}$ is bad, so is $\mu_{j-1}$, and hence, 
by applying Claim \ref{betano} to $j-1$ in place of $j$, we obtain that
$\mu_{j} > \mu_{j-1}$, from which we arrive at $|X_{j-1}| > |X_{j}|$, 
as required.  
\finbox

\medskip

Claim \ref{Xno} implies that 
$n \geq |X_{0}| > |X_{1}|> \cdots > |X_{h-1}| \geq 1$, from which $h\leq n$ follows.  
\finbox \finboxHere 
\medskip

\subsection{Computing a dec-min element in strongly polynomial time}
\label{strong.pol}

In order to compute a dec-min element of an M-convex set
$\odotZ{B}=\odotZ{B'}(p)$, 
our first task is to compute the smallest integer $\beta_{1}$ 
for which $\odotZ{B}$
has an element with largest component $\beta_{1}$.  
Theorem 4.1 of \cite{FM19partA} asserts that  
$\beta_{1} =\max \{ \lceil { p(X) / |X| }\rceil : \emptyset \not =X\subseteq S\}$.
By applying the ND-algorithm described in Section~\ref{pbmaximizing}, 
we can compute $\beta_{1}$ in strongly polynomial time.
Note that, by Theorem~\ref{Msteps}, the algorithm terminates after at most $n$ applications
of Subroutine \eqref{(ND.routine)}.

For any number  $\beta$, a vector is said to be {\bf $\beta$-covered} 
if each of its components is at most $\beta$.
An element $m$ of $\odotZ{B}$ is called a 
{\bf max-minimizer}
if its largest component is as small as possible.  
A max-minimizer element $m$ is said to be
{\bf pre-dec-min} in $\odotZ{B}$ 
if the number of its largest components is as small as possible.  
Obviously, a dec-min element is pre-dec-min, 
and a pre-dec-min element is a max-minimizer.

Given the value of $\beta_{1}$,  a $\beta_{1}$-covered element
 $m$ of $\odotZ{B}$
can easily be computed with a greedy-type algorithm as follows.
Since there is a $\beta_{1}$-covered member of $B$, the vector 
$(\beta_{1},\beta_{1},\dots ,\beta_{1})$ 
belongs to the so-called supermodular polyhedron 
$S'(p):=\{x:  \widetilde x(X)\geq p(X)$ for every $X\subseteq S\}$.  
Consider the elements of $S$ in an arbitrary order $\{s_{1},\dots ,s_{n}\}$.  
Let $m(s_{1}):=\min \{z:  (z,\beta_{1},\beta_{1},\dots ,\beta_{1})\in S'(p)\}$. 
In the general step, if the components 
$m(s_{1}),\dots ,m(s_{i-1})$ have already been determined, let
\begin{equation} 
m(s_{i}):= \min \{z: 
 (m(s_{1}),m(s_{2}),\dots ,m(s_{i-1}), z,\beta_{1},\beta_{1},\dots ,\beta_{1})\in S'(p)\}.  
\label{(m.def)} 
\end{equation}
This computation can be carried out 
by $n$ applications of a subroutine for a submodular function minimization.

Given a $\beta_{1}$-covered integral element $m$ of $B$, 
our next goal is to modify $m$ by a series of 1-tightening steps 
to obtain a pre-dec-min element of $\odotZ{B}$.  
To this end, we consider one by one 
those elements $t$ of $S$ for which 
$m(t) = \beta_{1}$, and check whether 
a 1-tightening step can be applied at $t$ with some element $s\in S$. 
That is, 
for each element $s\in S$ with $m(s) \leq \beta_{1}-2$, 
we check whether there is no $m$-tight $t\overline{s}$-set.  
If we find such an $s$, then $m':=m - \chi_{t}+\chi_{s}\in \odotZ{B}$ and $m'$ 
has one less coordinates with value $\beta_{1}$.
(Note that the largest component of $m'$ is also $\beta_{1}$ 
as $\beta_1$ was chosen to be the smallest upper bound).  
If no such element $s$ exists, then
\begin{equation}
 \hbox{ $\beta_{1}\geq m(s)\geq \beta_{1}-1$ holds for each $s\in T_m(t)$, }\ 
\label{(betaegy)} 
\end{equation}
where $T_m(t)$ is the smallest $m$-tight set containing $t$.  
While updating $m:=m'$ in the first case,
we iterate the above procedure for each element $t$ of $S$ with $m(t)=\beta_{1}$.

The element of $\odotZ{B}$ obtained from this series of modifications,
 which is denoted by $m$, has the property \eqref{(betaegy)} 
for all $t \in S$ with $m(t)=\beta_{1} $.
Theorem 4.2 of \cite{FM19partA} states 
that a $\beta_{1}$-covered element $m$ of $\odotZ{B}$ is pre-dec-min
precisely if $m(s)\geq \beta_{1}-1$ for each $s\in S_1(m)$, where
$S_{1}(m) =\cup \{T_{m}(t):  m(t)=\beta_{1}\}$.
This and \eqref{(betaegy)} imply that the final vector 
obtained by the above
procedure is a pre-dec-min element of $\odotZ{B}$ indeed.

Note that for a given pair $( t,s )$ 
of elements, 
deciding whether a 1-tightening step is applicable or not
can be done by a single submod-minimization, 
and hence 
deciding whether a given $t$ with $m(t)=\beta_{1}$ 
admits a 1-tightening step reducing $m(t)$
can be done by at most $n$ calls of a submod-minimizer.  
Altogether, the procedure above needs at most $n$ 1-tightening steps 
which can be carried out by $n\sp{2}$ calls of a submod-minimizer.

Recall that $T_{m}(t)$ denoted the unique smallest tight set containing $t$ 
when $p$ is (fully) supermodular.
But $T_{m}(t)$ can be described without explicitly referring to $p$ 
since an element $s\in S$ belongs to $T_{m}(t)$ precisely if 
$m':=m - \chi_{t} + \chi_{s}$ is in $B$, 
and this is computable 
by the subroutine \eqref{(st.step)}.  
In particular,  we can compute $S_{1}(m)$ by \eqref{(st.step)}.  
We use a short-hand notation $S_{1} := S_{1}(m)$.

Let $B_{1}$ denote the restriction of the base-polyhedron $B$ to 
$S_{1}$ and $B_{1}'$ the contraction of $B$ by $S_{1}$.
Theorem 4.6 of \cite{FM19partA}
states that, 
for $m_{1} \in {\bf Z}\sp{S_{1}}$ and 
$m_{1}' \in {\bf Z}\sp{S - S_{1}}$, 
$(m_{1},m_{1}')$ is a dec-min element of
$\odotZ{B}$ precisely if $m_{1}$ is a dec-min element of $\odotZ{B_1}$
and $m_{1}'$ is a dec-min element of $\odotZ{B_1'}$.
Let $m_{1}:=m\vert S_{1}$ for the pre-dec-min element $m$ constructed above.
Since $m_{1}$ is near-uniform 
on $S_{1}$, it is a dec-min element of $\odotZ{B_{1}}$.
Hence, if $m_{1}'$ is a dec-min element of $\odotZ{B_{1}'}$, 
then $(m_{1},m_{1}')$ is a dec-min element of $\odotZ{B}$.
Such a dec-min element $m_{1}'$ can be computed 
by applying iteratively the computation described above for computing $m_{1}$.
In this way we can compute a dec-min element of $\odotZ{B}$.

It is worth mentioning that 
the restriction of the pre-dec-min element $m$ to 
$S_{1} = S_{1}(m)$
is the same as 
the restriction of the dec-min element (found by the above algorithm) to $S_1$. 
Therefore this $S_1$ is the first member of 
the canonical partition belonging to $\odotZ{B}$.

The above algorithm computes a dec-min element of $\odotZ{B}$ in strongly polynomial time.
The subroutine \eqref{(member)} 
is called  only once at the beginning of the algorithm,
and the running time of the algorithm 
is governed by the number of calls of 
\eqref{(st.step)} and \eqref{(ND.routine)}. 
First assume that we are given a fully supermodular function $p$ to describe $B$.
To determine $m_1$ and $S_{1}$,  we need
(i)
$n$ calls of \eqref{(ND.routine)} to compute the value $\beta_1$ by the ND-algorithm, 
where one call of \eqref{(ND.routine)} requires
a single application of a submod-minimizer,
(ii) $n$ applications of a submod-minimizer 
 to compute a $\beta_1$-covered element by the greedy-type algorithm, and  
(iii) $n\sp{2}$ calls of \eqref{(st.step)}
 to compute a pre-dec-min element,
 where one call of \eqref{(st.step)} requires
 a single application of a submod-minimizer. 
Therefore, we can determine $m_1$ and $S_{1}$
with $n+n+ n\sp{2} = O(n\sp{2})$ applications of a submod-minimizer.
We repeat the above procedure on $S - S_{1}$, $S - (S_{1} \cup S_{2})$, and so on.
Hence the above algorithm finds
a dec-min element and the canonical partition
with $O(n\sp{3})$ applications of a submod-minimizer,
provided that we are given a fully supermodular function $p$ to describe 
the base-polyhedron $B$.
Even when $B$ is given in terms of an intersecting or crossing supermodular function,
this algorithm is strongly polynomial.
Namely, the total number of submod-minimizer calls
is $O(n\sp{4})$ if $p$ is intersecting supermodular
and $O(n\sp{5})$ if $p$ is crossing supermodular
(see Remark \ref{RMinsercross}).
Finally we recall Remark \ref{RMsubmmin} for the complexity of a submod-minimizer.

\section{Applications} 
\label{alk}

\subsection{Background}

There are two major sources of applicability of the 
structural results on decreasing minimization on an M-convex set.
One of them relies on the fact that the class of
integral base-polyhedra is closed under several operations.  
For example, a face of a base-polyhedron is also a base-polyhedron, and so
is the intersection of an integral box with a base-polyhedron $B$.
Also, the sum of integral base-polyhedra $B_{1},\dots ,B_{k}$ is a
base-polyhedron $B$ which has, in addition, the integer decomposition
property meaning that any integral element of $B$ can be obtained as
the sum of $k$ integral elements by taking one from each $B_{i}$.  
This latter property implies that the sum of M-convex sets is M-convex.  
We also mention the important operation of taking an aggregate of a
base-polyhedron, to be introduced below in Section \ref{matroidalk}.

The other source of applicability is based on the fact that not only
fully super- or submodular functions can define base-polyhedra but
some weaker functions as well.  
For example, if $p$ is an integer-valued crossing 
(in particular, intersecting) supermodular
function with finite $p(S)$, then $B=B'(p)$ is a (possibly empty)
integral base-polyhedron (and $\odotZ{B}$ is an M-convex set).  
This fact will be exploited in solving dec-min orientation problems when
both degree-constraints and edge-connectivity requirements must be fulfilled.  
In some cases even weaker set-functions can define base-polyhedra.  
This is why we can solve dec-min problems concerning
edge- and node-connectivity augmentations of digraphs.

\subsection{Applications to matroids} 
\label{matroidalk}

Levin and Onn \cite{Levin-Onn} solved algorithmically the following problem:
\textit{Find $k$ bases of a matroid $M$ on a ground-set $S$ such
that the sum of their characteristic vectors be decreasingly minimal.}
Their approach, however, does not seem to work in the following natural extension.  
Suppose we are given $k$ matroids $M_{1},\dots ,M_{k}$
on a common ground-set $S$, and our goal is to find a basis $B_{i}$ of
each matroid $M_{i}$ in such a way that the vector  
$\sum [\chi_{B_{i}} : i=1,\dots ,k]$ 
is decreasingly minimal.  
Let 
$B_{\Sigma}$ 
denote the sum of the base-polyhedra of the $k$ matroids.  
By a theorem of 
Edmonds \cite{Edmonds70},
 the integral elements of 
$B_{\Sigma}$ 
are exactly the vectors of form 
$\sum [\chi_{B_{i}} :  i=1,\dots ,k]$ 
where $B_{i}$ is a basis of $M_{i}$.  
Therefore the problem is to find a dec-min element of 
$\odotZ{B_{\Sigma}}$.  
This can be done
 by the basic algorithm described in Section \ref{SCbasicalg}.
Let us see how the requested subroutines are available in this special case.  
The algorithm starts with an arbitrary member $m$ of 
$\odotZ{B_{\Sigma}}$
which is obtained by taking a basis $B_{i}$ from each matroid $M_{i}$, 
and these bases define $m:=\sum_{i} \chi_{B_{i}}$.

To realize Subroutine \eqref{(st.step)}, we mentioned that it suffices
to realize Subroutine \eqref{(st.step2)}, which requires for a given
integral vector $m'$ with $\widetilde m'(S)= \sum_{i}r_{i}(S)$ to decide
whether $m'$ is in 
$\odotZ{B_{\Sigma}}$
 or not.  
But this can simply be done by Edmonds' matroid intersection algorithm
\cite{Edmonds79} (see also Section 13.1.2 in \cite{Frank-book}.
Namely, let
$S_{1},\dots ,S_{k}$ be disjoint copies of $S$ and $M_{i}'$ an isomorphic
copy of $M_{i}$ on $S_{i}$.  
Let $N_{1}$ be the direct sum of matroids
$M_{i}'$ on ground-set $S':=S_{1}\cup \cdots \cup S_{k}$.  
Let $N_{2}$ be a partition matroid on $S'$ in which a subset $Z$ is a basis 
if it contains exactly $m'(s)$ members of the $k$ copies of $s$ 
for each $s\in S$.  
Then $m'$ is in 
$\odotZ{B_{\Sigma}}$ precisely if $N_{1}$ and
$N_{2}$ have a common basis.

In conclusion, with the help of Edmonds' matroid intersection algorithm, 
Subroutine \eqref{(st.step)} is available, and hence the
basic algorithm can be applied.

\medskip \medskip

Another natural problem concerns a single matroid $M$ on a ground-set $T$.  
Suppose we are given a partition ${\cal P}= \{T_{1},\dots ,T_{n}\}$
of $T$ and we consider the intersection vector 
$(\vert Z\cap T_{1}\vert,\dots ,\vert Z\cap T_{n}\vert )$ 
assigned to a basis $Z$ of $M$.  
The problem is to find a basis for which the intersection vector is decreasingly minimal.

To solve this problem, we recall an important construction of
base-polyhedra, called the aggregate.  
Let $T$ be a ground-set and $B_{T}$ an integral base-polyhedron in ${\bf R}\sp T$.  
Let ${\cal P}=\{T_{1},\dots ,T_{n}\}$ be a partition of $T$ into non-empty subsets
and let $S=\{s_{1},\dots ,s_{n}\}$ be a set 
whose elements correspond to the members of ${\cal P}$.  
The aggregate $B_{S}$ of $B_{T}$ is defined as follows.
\begin{equation} 
\hbox{ $B_{S}:= \{ (y_{1},\dots ,y_{n}):$ there is an $x\in B_{T}$ with
   $y_{i}=\widetilde x(T_{i}) \ (i=1,\dots ,n)\}$.   }\ 
\label{(aggre)} 
\end{equation}
\noindent
A basic theorem concerning base-polyhedra 
(see, for example, Theorems 14.2.12 and 14.2.13 in book \cite{Frank-book})
states that $B_{S}$ is a
base-polyhedron, moreover, for each integral member $(y_{1},\dots ,y_{n})$
of $B_{S}$, the vector $x$ in \eqref{(aggre)} can be chosen integer-valued.  
In other words,
\begin{equation} 
\hbox{ $\odotZ{B_{S}}:= \{ (y_{1},\dots ,y_{n}):$ there is an $x\in
\odotZ{B_{T}}$ with $y_{i}=\widetilde x(T_{i}) \ (i=1,\dots ,n)\}$.   }\
\label{(aggre.odot)} 
\end{equation}
\noindent
We call $\odotZ{B_{S}}$ the {\bf aggregate} of $\odotZ{B_{T}}$.

Returning to our matroid problem, let $B_{T}$ denote the base-polyhedron
of matroid $M$.  Then the problem is nothing but finding a dec-min
element of $\odotZ{B_{S}}$.

We can apply the basic algorithm (concerning M-convex sets) for this
special case since the requested subroutines are available through
standard matroid algorithms.  Namely, Subroutine \eqref{(member)} is
available since for any basis $Z$ of $M$, the intersection vector
assigned to $Z$ is nothing but an element of $\odotZ{B_{S}}$.

To realize Subroutine \eqref{(st.step)}, we mentioned that it suffices
to realize Subroutine \eqref{(st.step2)}.  
Suppose we are given a vector $y\in {\bf Z}_{+}\sp{S}$ 
(Here $y$ stands for $m'$ in \eqref{(st.step2)}).  
Suppose that $\widetilde y(S)=r(T)$
 (where $r$ is the rank-function of matroid $M$) and that 
$y(s_{i})\leq \vert T_{i}\vert$ for $i=1,\dots ,n$.

Let $G=(S,T;E)$ denote a bipartite graph where $E=\{ts_{i}:  t\in T_{i},
i=1,\dots ,n\}$.  By this definition, the degree of every node in $T$
is 1 and hence the elements of $E$ correspond to the elements of $M$.
Let $M_{1}$ be the matroid on $E$ corresponding to $M$ (on $T$).  Let
$M_{2}$ be a partition matroid on $E$ in which a set $F\subseteq E$ is a
basis if $d_{F}(s_{i})=y(s_{i})$,
where $d_{F}(s_{i})$ denotes the number of edges in $F$ for which $s_{i}$ is an end-node.
By this construction, the vector $y$ is
in $\odotZ{B_{S}}$ precisely if the two matroids $M_{1}$ and $M_{2}$ have a
common basis.  This problem is again tractable 
by Edmonds' matroid intersection algorithm.

\medskip 
As a special case, we can find a spanning tree of a
(connected) directed graph for which its in-degree-vector is
decreasingly minimal.  Since the family of unions of $k$ disjoint
bases of a matroid forms also a matroid, we can also compute $k$
edge-disjoint spanning trees in a digraph whose union has a
decreasingly minimal in-degree vector.

Another special case is when we want to find a spanning tree of a
connected bipartite graph $G=(S,T;E)$ whose in-degree vector
restricted to $S$ is decreasingly minimal.

\medskip

\subsection{Applications to flows} 
\label{SCflows}

\subsubsection{A base-polyhedron associated with net-in-flows}
\label{SCflowbase}

Let $D=(V,A)$ be a digraph endowed with integer-valued bounding functions 
$f:A\rightarrow {\bf Z}\cup \{-\infty \}$ and
$g:A\rightarrow {\bf Z}\cup \{+\infty \}$ for which $f\leq g$.  
We call a vector (or function) $z$ on $A$ {\bf feasible} 
if $f\leq z\leq g$.  
The {\bf net-in-flow} $\Psi_z$ of $z$ is a vector on $V$ and defined by 
$\Psi_z(v) = \varrho_z(v) -\delta_z(v)$, 
where 
$\varrho_z(v):=\sum [z(uv):  uv\in A]$ and $\delta_z(v):=\sum [z(vu):  uv\in A]$.
If $m$ is the net-in-flow of a vector $z$, 
then we also say that $z$ is an {\bf $m$-flow}.

A variation of Hoffman's classic theorem on feasible circulations
\cite{Hoffman60} is as follows.

\begin{lemma} 
\label{LMnetinflowbase}
An integral vector $m:V\rightarrow {\bf Z} $ is the net-in-flow
of an integral feasible vector 
(or in other words, there is an integer-valued feasible $m$-flow) 
if and only if 
$\widetilde m(V)=0$ and
\begin{equation} 
\varrho_f(Z) - \delta_g(Z) \leq \widetilde m(Z)
 \ \ \hbox{\rm holds whenever}\ \ \ Z\subseteq V, 
\label{(feasmflow)} 
\end{equation}
where \ 
$\varrho_f(Z):=\sum [f(a):  a\in A \hbox{\  \rm and $a$ enters $Z$} ]$ 
\ and \ 
$\delta_g(Z):=\sum [g(a):  a\in A \hbox{\  \rm and $a$ leaves $Z$}]$.
\finbox 
\end{lemma}

Define a set-function $p_{fg}$ on $V$ by 
\[
  p_{fg}(Z):=\varrho_f(Z)-\delta_g(Z).
\] 
Then $p_{fg}$ is (fully) supermodular 
(see, e.g.  Proposition 1.2.3 in \cite{Frank-book}).  
Consider the base-polyhedron
$B_{fg}:=B'(p_{fg})$ and the M-convex set $\odotZ{B_{fg}}$.  
By Lemma \ref{LMnetinflowbase} 
the M-convex set $\odotZ{B_{fg}}$ consists
exactly of the net-in-flow integral vectors $m$.

By the algorithm described in Section \ref{SCalgo2}, we can compute a
decreasingly minimal element of $\odotZ{B_{fg}}$ in strongly polynomial time.  
By relying on a variant of the strongly polynomial push-relabel
algorithm described in book \cite{Frank-book} (see, Section 6.1.3),
one can check whether or not \eqref{(feasmflow)} holds.  If it does
not, then this variant can compute a set most violating
\eqref{(feasmflow)} 
(that is, a maximizer of $\varrho_f(Z) - \delta_g(Z) - \widetilde m(Z)$), 
while if \eqref{(feasmflow)} does hold, 
then the push-relabel algorithm computes an integral valued feasible $m$-flow. 
Therefore the requested oracles in the general algorithm 
for computing a dec-min element are available through a
network flow algorithm, and we do not have to rely on a
general-purpose submodular function minimizing oracle.

For the sake of an application of this algorithm to capacitated
dec-min orientations in Section \ref{capori}, we remark that the
algorithm can also be used to compute a dec-min element of the
M-convex set obtained from $\odotZ{B_{fg}}$ 
by translating it with a given integral vector.

\subsubsection{Discrete version of Megiddo's flow problem}

Megiddo \cite{Megiddo74}, \cite{Megiddo77} considered the following problem.  
Let $D=(V,A)$ be a digraph endowed with a non-negative
capacity function 
$g:A\rightarrow {\bf R}_{+}$.  
Let $S$ and $T$ be two disjoint non-empty subsets of $V$.  
Megiddo described an algorithm to compute a feasible flow from $S$ to $T$ 
with maximum flow amount $M$
for which the net-in-flow vector restricted on $S$ 
is (in our terms) increasingly maximal.  
Here a feasible flow is a vector $x$ on $A$ 
for which $\Psi_{x}(v) \leq 0$ for $v\in S$, 
$\Psi_{x}(v)\geq 0$ for $v\in T$,
and $\Psi_{x}(v)=0$ for $v\in V-(S\cup T)$.  
The flow amount $x$ is $\sum [\Psi_{x}(t):t\in T]$.

We emphasize that Megiddo solved the continuous (fractional) case and
did not consider the corresponding discrete (or integer-valued) flow problem.  
To our knowledge, this natural optimization problem has not been investigated so far.

To provide a solution, suppose that $g$ is integer-valued.  
Let $f\equiv 0$ and consider the net-in-flow vectors belonging to feasible vectors.  
These form a base-polyhedron $B_{1}$ in ${\bf R}\sp V$.  
Let $B_{2}$ denote the 
base-polyhedron
obtained from $B_{1}$ by intersecting
it with the box defined by 
$z(v)\leq 0$ for $v\in S$, $z(v)\geq 0$ for $v\in T$ 
and $z(v)=0$ for $v\in V-(S\cup T)$.   

The restriction of $B_{2}$ to $S$ is a g-polymatroid $Q$ in ${\bf R}\sp{S}$.  
And finally, we can consider the face of $Q$ defined by $\widetilde z(S)=-M$.  
This is a base-polyhedron $B_3$ in ${\bf R}\sp{S}$, 
and the discrete version of Megiddo's flow problem is equivalent 
to finding an inc-max element of $\odotZ{B_3}$. 
(Recall that an element of an M-convex set is dec-min precisely if it is inc-max.)

It can be shown that in this case again 
the general submodular function minimizing subroutine used in the algorithm 
to find a dec-min element of an M-convex set can be replaced by a max-flow min-cut algorithm.

A recent paper \cite{FM19partC} addresses a more general problem 
to find an integral feasible flow
that is dec-min on an arbitrarily specified edge set.

\subsection{Further applications}

\subsubsection{Root-vectors of arborescences}

A graph-example comes from packing arborescences.  
Let $D=(V,A)$ be a digraph and $k>0$ an integer.  
We say that a non-negative integral vector 
$m:V\rightarrow {\bf Z}_{+}$
is a {\bf root-vector} if
there are $k$ edge-disjoint spanning arborescences such that each node
$v\in V$ is the root of $m(v)$ arborescences.  
Edmonds \cite{Edmonds73} classic result on disjoint arborescences implies that
$m$ is a root-vector 
if and only if 
$\widetilde m(V)=k$ and
$\widetilde m(X)\geq k-\varrho(X)$ holds for every subset $X$ with
$\emptyset \subset X\subset V$.  
Define set-function $p$ by $p(X):  = k-\varrho(X)$ 
if $\emptyset \subset X\subseteq V$ and $p(\emptyset):=0$.  
Then $p$ is intersecting supermodular, so $B'(p)$ is an integral base-polyhedron.  
The intersection $B$ of $B'(p)$ with the
non-negative orthant is also a base-polyhedron, and the theorem of
Edmonds is equivalent to stating that a vector $m$ is a root-vector 
if and only if 
$m$ is in $\odotZ{B}$.

Therefore the general results on base-polyhedra can be specialized to
obtain $k$ disjoint spanning arborescences whose root-vector is
decreasingly minimal.

\subsubsection{Connectivity augmentations}

Let $D=(V,A)$ be a directed graph and $k>0$ an integer.  
We are interested in finding a so-called augmenting digraph 
$H=(V,F)$ of $\gamma $ arcs for which $D+H$ is $k$-edge-connected or
$k$-node-connected.  
In both cases, the in-degree vectors of the
augmenting digraphs are the integral elements of an integral
base-polyhedron \cite{FrankJ23}, \cite{FrankJ31}.  
Obviously, the in-degree vectors of the augmented digraphs 
are the integral elements of an integral base-polyhedron.

Again, our results on general base-polyhedra can be specialized to
find an augmenting digraph whose in-degree vector is decreasingly minimal.

\section{Orientations of graphs} 
\label{SCori1}

The literature is quite rich in graph orientation problems 
where the task is to orient the edges of an undirected graph 
such that the resulting digraph meet some expected properties.  
For an overview of graph orientation problems, 
see, for example, Chapter~61 of the book of Schrijver \cite{Schrijverbook} 
or Chapter~9 of the book of Frank \cite{Frank-book}.  
This latter one illuminates the deep connection 
between orientation problems and submodular optimization.
Although quite general and efficient tools are exhibited in these books 
to manage graph orientation problems, 
they do not say anything 
about the orientation problem investigated first 
in a recent paper of Borradaile et al.~\cite{BIMOWZ},
which is the decreasingly minimal orientation problem in our terms. 
It was actually their paper that triggered 
the whole research behind our present work.  
In this section, we exhibit how the theoretical background 
developed in \cite{FM19partA} and the algorithms of 
Section \ref{SCalgo2} can be used to obtain
major extensions of results in \cite{BIMOWZ}.

\medskip

Let $G=(V,E)$ be an undirected graph.  
For $X\subseteq V$, let $i_{G}(X)$
denote the number of edges induced by $X$
while $e_{G}(X)$
is the number of edges with at least one end-node in $X$.  
Then $i_{G}$ is supermodular, $e_{G}$ is submodular, and they are
complementary functions, that is, 
$i_{G}(X) = e_{G}(V)- e_{G}(V-X)$.  
Let $B_{G}:= B(e_{G}) = B'(i_{G})$ 
denote the base-polyhedron defined by $e_{G}$ or $i_{G}$.

We say that a function $m:V\rightarrow {\bf Z} $ is 
the {\bf in-degree vector} of an orientation $D$ of $G$ 
if $\varrho_{D}(v)=m(v)$ for each node $v\in V$.  
An in-degree vector $m$ obviously meets the equality
$\widetilde m(V)=\vert E\vert $. 
The following basic result, sometimes called the Orientation lemma, 
is due to Hakimi \cite{Hakimi}.

\begin{lemma}[Orientation lemma] \label{LMorientationlemma}
Let $G=(V,E)$ be an undirected graph and 
$m: V\rightarrow {\bf Z}$ an integral vector for which 
$\widetilde m(V)=\vert E\vert$. 
Then $G$ has an orientation with in-degree vector $m$ 
if and only if
\begin{equation} 
\widetilde m(X) \leq e_{G}(X)\ \ \hbox{\rm for every subset \ $X\subseteq V$,}\ 
\label{(eG)} 
\end{equation}
which is equivalent to
\begin{equation} 
\widetilde m(X) \geq i_{G}(X) \ \  \hbox{\rm for every subset \ $X\subseteq V$.}\  
\label{(iG)} 
\end{equation}
\vspace{-1\baselineskip}\\
\finbox
\end{lemma}

This immediately implies the following claim.

\begin{claim} \label{oribase} 
The in-degree vectors of orientations of $G$
are precisely the integral elements of base-polyhedron 
$B_{G}$\ $(=B(e_{G}) = B'(i_{G}))$, 
that is, the set of in-degree vectors of orientations of $G$ 
is the M-convex set $\odotZ{B_{G}}$.  
\finbox
\end{claim}

The proof of Lemma \ref{LMorientationlemma} is algorithmic 
(see, e.g.,  Theorem 2.3.2 of \cite{Frank-book})
and the orientation corresponding to a given $m$ can be constructed easily.

\subsection{Decreasingly minimal orientations}
\label{SCdecminori}

Due to Claim \ref{oribase}, we can apply the results on
dec-min elements to the special base-polyhedron $B_{G}$.  
Borradaile et al.~\cite{BIMOWZ} called an orientation of $G$
egalitarian if its in-degree vector is decreasingly minimal 
but we prefer the term {\bf dec-min orientation} 
since an orientation with an increasingly maximal in-degree vector 
also has an intuitive
egalitarian feeling.  Such an orientation is called {\bf inc-max}.
Theorem~\ref{THequidotone} immediately implies the following.

\begin{corollary}  \label{dm-im} 
An orientation of $G$ is dec-min 
if and only if 
it is inc-max.  
\finbox 
\end{corollary}

Note that the term dec-min orientation is asymmetric in the sense that
it refers to in-degree vectors.  
One could also aspire for finding an orientation 
whose out-degree vector is decreasingly minimal.  
But this problem is clearly equivalent to the in-degree version and hence in
the present work we do not consider out-degree vectors, 
apart from a single exception in Section \ref{inout}.

By Theorem~\ref{THequidotone},  an element $m$ of $\odotZ{B_{G}}$ 
is decreasingly minimal 
if and only if 
there is no 1-tightening step for $m$.  
What is the meaning of a 1-tightening step in terms of orientations?

\begin{claim} \label{reori} 
Let $D$ be an orientation of $G$ with in-degree vector $m$.  
Let $t$ and $s$ be nodes of $G$.  
The vector $m':=m + \chi_{s} - \chi_{t}$ is in $B_{G}$ 
if and only if 
$D$ admits a dipath from $s$ to $t$.  
\end{claim}  

\Proof 
$m'\in B_{G}$ holds precisely if there is no $t\overline{s}$-set $X$
which is tight with respect to $i_{G}$, 
that is, $\widetilde m(X) = i_{G}(X)$.  
Since 
$\varrho(Y) + i_{G}(Y) = \sum [\varrho(v):v\in Y]  = \widetilde m(Y)$
holds for any set $Y\subseteq V$, 
the tightness of $X$ is equivalent to requiring that 
$\varrho(X)=0$.  
Therefore $m'\in B_{G}$ 
if and only if 
$\varrho(Y) > 0$ holds for every $t\overline{s}$-set $Y$, 
which is equivalent to the existence of a dipath of $D$ from $s$ to $t$.  
\finbox 
\medskip

Recall that a 1-tightening step at a member $m$ of $B_{G}$ consists of
replacing $m$ by $m'$ provided that $m(s)\geq m(t)+2$ and $m'\in B_{G}$.
By Claim \ref{reori}, a 1-tightening step at a given orientation of
$G$ corresponds to reorienting an arbitrary dipath 
from a node $s$ to node $t$  for which 
$\varrho(s)\geq \varrho(t)+2$.  
Therefore, Theorem~\ref{THequidotone}
immediately implies 
the following basic theorem of Borradaile et al.~\cite{BIMOWZ}.

\begin{theorem}[Borradaile et al.~\cite{BIMOWZ}] \label{pathrev} 
An orientation $D$ of a graph $G=(V,E)$ is decreasingly minimal 
if and only if 
no dipath exists from a node $s$ to a node $t$ 
for which $\varrho(t)\geq \varrho(s)+2$.  
\finbox 
\end{theorem}  

Note that this theorem also implies Corollary \ref{dm-im}.  
It immediately gives rise to an algorithm for finding a dec-min orientation.  
Namely, we start with an arbitrary orientation of $G$.
We call a dipath {\bf feasible} if $\varrho(t)\geq \varrho(s)+2$
holds for its starting node $s$ and end-node $t$.  
The algorithm consists of reversing feasible dipaths as long as possible.  
Since the sum of the squares of in-degrees always drops 
when a feasible dipath is reversed, 
and originally this sum is at most $\vert E\vert \sp{2}$,
the dipath-reversing procedure terminates after at most 
$\vert E\vert \sp{2}$ reversals.  
By Theorem~\ref{pathrev}, when no more feasible dipath exists, 
the current orientation is dec-min.  
The basic algorithm concerning general base-polyhedra in 
Section~\ref{SCbasicalg}
is nothing but an extension of the algorithm of Borradaile et al.

It should be noted that they suggested to choose at every step the
current feasible dipath in such a way that the in-degree of its
end-node $t$ is as high as possible, and they proved that the
algorithm in this case terminates after 
at most $\vert E\vert \vert V\vert $ dipath reversals.

Note that we obtained Corollary \ref{dm-im} as a special case of a
result on M-convex sets but it is also a direct consequence 
of Theorem~\ref{pathrev}.

\medskip

\subsection{Capacitated orientation} 
\label{capori}

Consider the following capacitated version 
of the basic dec-min orientation problem of Borradaile et al.~\cite{BIMOWZ}.  
Suppose that a positive integer $\ell(e)$ 
is assigned to each edge $e$ of $G$.  
Denote by $G\sp +$ the graph arising from $G$ 
by replacing each edge $e$ of $G$ with $\ell(e)$ parallel edges.  
Our goal is to find a dec-min orientation of $G\sp +$.  
In this case, an orientation of $G\sp +$ is described 
by telling that, among the $\ell(e)$ parallel
edges connecting the end-nodes $u$ and $v$ of $e$ 
how many are oriented toward $v$ 
(implying that the rest of the $\ell(e)$ edges are oriented toward $u$).  
In principle, this problem can be solved by applying the algorithm 
described above to $G\sp +$, and this algorithm
is satisfactory when $\ell$ is small in the sense that its largest
value can be bounded by a power of $\vert E\vert $. 
The difficulty in the general case is that the algorithm will be polynomial 
only in the number of edges of $G\sp +$, 
that is, in $\widetilde \ell(E)$, 
and hence this algorithm is not polynomial in $\vert E\vert $.

We show how the algorithm in 
Section~\ref{SCflowbase}
can be used to solve
the decreasingly minimal orientation problem in the capacitated case
in strongly polynomial time.  
To this end, let $D=(V,A)$ be an arbitrary orientation of $G$ 
serving as a reference orientation.
Define a capacity function $g$ on $A$ by 
$g(\vec{e}):=\ell(e)$,
where $\vec{e}$ denotes the arc of $D$ obtained by orienting $e$.

We associate an orientation of $G\sp +$ with an integral vector 
$z:A\rightarrow {\bf Z}_{+}$
with $z\leq g$ as follows.  For an arc $uv$ of $D$, 
orient $z(uv)$ parallel copies of $e=uv\in E$ toward $v$ and
$g(uv)-z(uv)$ parallel copies toward $u$.  
Then the in-degree of a node $v$ is 
$m_z(v):= \varrho_z(v) + \delta_{g-z}(v) 
  = \varrho_z(v) - \delta_z(v) +\delta_g(v)$.
Therefore our goal is to find an integral vector $z$ on $A$ 
for which $0\leq z\leq g$ and the vector $m_z$ on $V$ is dec-min.  
Consider the set of net-in-flow vectors
$\{(\Psi_z(v):  v\in V) :  0\leq z\leq g\}$.  
In Section~\ref{SCflowbase},
we proved that this is a base-polyhedron $B_{1}$.  
Therefore the set of vectors $(m_z(v):v\in V)$ is also a base-polyhedron $B$ arising from
$B_{1}$ by translating $B_{1}$ with the vector $(\delta_g(v):v\in V)$.

As remarked at the end of Section~\ref{SCflowbase},
a dec-min element of
$\odotZ{B}$ can be computed in strongly polynomial time 
by a variant (Section 6.1.3 in \cite{Frank-book}) of the 
push-relabel subroutine for network flows 
(and not using a general-purpose submodular function minimizer).

\subsection{Canonical chain and essential value-sequence for orientations}
\label{SCoricano}

In Section \ref{SCcompcanchain}, we described Algorithm \ref{algo1} for an arbitrary
M-convex set $\odotZ{B}$ that computes, from a given dec-min element
$m$ of $\odotZ{B}$, the canonical chain and essential value-sequence
belonging to $\odotZ{B}$.  
That algorithm needed an oracle for computing 
the smallest $m$-tight set $T_{m}(u)$ containing $u$.  
Here we show how this general algorithm can be turned into 
a pure graph-algorithm in the special case of dec-min orientations.

To this end, consider the special M-convex set, denoted by $\odotZ{B_{G}}$, 
consisting of the in-degree vectors of the orientations 
of an undirected graph $G=(V,E)$.  
By the Orientation lemma, 
$B_{G}=B'(i_{G})$ where $i_{G}(X)$ denotes the number of edges induced by $X$.  
Recall that $i_{G}$ is a fully supermodular function.  
For an orientation $D$ of $G$ with in-degree vector $m$, the smallest
$m$-tight set $T_{m}(t)$ (with respect to $i_{G}$) containing a node $t$
will be denoted by $T_{D}(t)$.

\begin{claim} \label{mtight} 
Let $D$ be an arbitrary orientation of $G$ with in-degree vector $m$.  
{\rm (A)}\  A set $X\subseteq V$ is $m$-tight (with respect to $i_{G}$) 
if and only if $\varrho_{D}(X)=0$. 
{\rm (B)} \ The smallest $m$-tight set $T_{D}(t)$ containing a node $t$ 
is the set of nodes from which $t$ is reachable in $D$.
\end{claim}

\Proof We have
\[
 \varrho_{D}(X) + i_{G}(X) 
 = \sum [\varrho_{D}(v):  v\in X] = \widetilde m(X) \geq i_{G}(X) ,
\]
from which $X$ is $m$-tight 
(that is, $\widetilde m(X) =i_{G}(X)$) precisely if $\varrho_{D}(X)=0$, 
and Part (A) follows.
Therefore the smallest $m$-tight set $T_{D}(t)$ containing $t$ is the smallest set 
containing $t$ with in-degree 0, and hence $T_{D}(t)$ is
indeed the set of nodes from which $t$ is reachable in $D$, as stated
in Part (B).  
\finbox \medskip

By Claim \ref{mtight}, $T_{D}(t)$ is easily computable, and hence
Algorithm \ref{algo1} for general M-convex sets can easily be
specialized to graph orientations.  
By applying Theorem~\ref{THchardecmin}
to $p:=i_{G}$ and recalling from Claim \ref{mtight} that $C_{i}$ is
$m$-tight, in the present case, precisely if $\varrho_{D}(C_{i})=0$, 
we obtain the following.

\begin{theorem} \label{decminori} 
An orientation of $D$ of $G$ is dec-min if
and only of $\varrho_{D}(C_{i})=0$ for each member $C_{i}$ of the canonical
chain and $\beta_{i} -1\leq \varrho_{D}(v)\leq \beta_{i}$ holds for every
node $v\in S_{i}$ $(i=1,\dots ,q)$.  
\finbox 
\end{theorem}

We remark that the members of the canonical partition computed by our
algorithm for $\odotZ{B_{G}}$ is exactly the non-empty members of the
so-called density decomposition of $G$ introduced by Borradaile et al.~\cite{BMW}.
We also remark that by combining the approach of the present section
with Section \ref{capori}, the canonical partition and the essential
value-sequence can be computed in strongly polynomial time even in the
capacitated orientation problem.

\subsection{Cheapest dec-min orientations}
\label{cheapdecmin}

It is indicated in \cite{FM19partA} that, 
in decreasing minimization on an M-convex set in general, 
we can construct an algorithm to compute a cheapest dec-min element 
with respect to a given (linear) cost-function on the ground-set. 
In the special case of dec-min orientations, 
this means that if $c$ is a cost-function on the node-set of $G=(V,E)$, 
then we have an algorithm to compute a dec-min orientation of $G$ for which 
$\sum [c(v)\varrho(v):v\in V]$ is minimum.

But the question remains:  what happens if, instead of a cost-function
on the node-set, we have a cost-function $c$ on $\vec E_{2}$, where
$\vec E_{2}$ arises from $E$ by replacing each element $e=uv$ ($=vu$) of
$E$ by two oppositely oriented arcs $uv$ and $vu$, and we are
interested in finding a cheapest orientation with specified properties?
  (As an orientation of $e$ consists of replacing $e$ by
one of the two arcs $uv$ and $vu$ and the cost of its orientation is,
accordingly, $c(uv)$ or $c(vu)$.  
Therefore we can actually assume that $\min \{c(uv),c(vu)\}=0.$)

It is important to remark that the standard minimum cost 
in-degree specified orientation problem can be easily reduced, 
with a straightforward technique, to a minimum cost flow problem in a digraph
with small integral capacities 
(see, e.g.  Section 3.6.1 in book \cite{Frank-book}).  
The same reduction works for min-cost in-degree constrained orientations, as well.  
Note that already the min-cost flow algorithm of 
Ford and Fulkerson \cite{Ford-Fulkerson} is strongly polynomial 
when the capacities are small integers 
(that is, we do not need here the significantly more sophisticated min-cost flow algorithm
of Tardos \cite{Tardos1} which is strongly polynomial for arbitrary capacities.)  
Actually, we shall need a version of this minimum cost orientation problem when
some of the edges are already oriented, and this slight extension is
also tractable by network flows.

Theorem~\ref{decminori} implies that 
the problem of finding a cheapest dec-min orientation is equivalent to 
finding a cheapest 
in-degree constrained orientation
by orienting edges connecting $C_{i}$ and $V-C_{i}$ toward $V-C_{i}$
($i=1,\dots ,q$).  
Here the in-degree constraints are given by 
$\beta_{i}-1\leq \varrho_{D}(v)\leq \beta_{i}$ for $v\in S_{i}$ ($i=1,\dots ,q$).

Note that Harada et al.~\cite{HOSY07} provided a direct algorithm for
 the minimum cost version of the so-called semi-matching problem, 
which problem includes the minimum cost dec-min orientation problem. 
For this link, see Section \ref{SCsemimatching}.

We remark that by combining the approach of the present section with Section \ref{capori}, 
a cheapest dec-min orientation can be computed in strongly polynomial time 
even in the capacitated orientation problem.
In this case, however, one needs a strongly polynomial subroutine to
compute a minimum cost feasible circulation.  The first such algorithm
is due to Tardos \cite{Tardos1}.

\subsection{Orientation with dec-min in-degree vector and dec-min out-degree vector} 
\label{inout}

We mentioned that dec-min and inc-max orientations always concern in-degree vectors.  
As an example to demonstrate the advantage of the general base-polyhedral view, 
we outline here one exception when
in-degree vectors and out-degree vectors play a symmetric role.  
The problem is to characterize undirected graphs admitting an orientation
which is both dec-min with respect to its in-degree vector and dec-min
with respect to its out-degree vector.

For the present purposes, we let $d_{G}$ denote the degree vector of $G$, 
that is, $d_{G}(v)$ is the number of edges incident to $v\in V$.
(This notation differs from the standard set-function meaning of $d_{G}$.)

Let $B_{\rm in}$ denote the convex hull of the in-degree vectors of
orientations of $G$, and $B_{\rm out}$ the convex hull of out-degree
vectors of orientations of $G$. 
 (Earlier $B_{\rm in}$ was denoted by
$B_{G}$ but now we have to deal with both out-degrees and in-degrees.)
As before, $\odotZ{B_{\rm in}}$ 
is the set of in-degree vectors of
orientations of $G$, and $\odotZ{B_{\rm out}}$ is the set of
out-degree vectors of orientations of $G$. 
Let 
$\odotZ{B_{\rm in}\sp{\bullet}}$ 
denote the set of dec-min in-degree vectors of orientations
of $G$, and 
$\odotZ{B_{\rm out}\sp{\bullet}}$ 
the set 
of dec-min out-degree vectors of orientations of $G$.  
By Theorem~\ref{THmatroid-eltolt},
both $\odotZ{B_{\rm in}\sp{\bullet}}$ and
$\odotZ{B_{\rm out}\sp{\bullet}}$ are matroidal M-convex sets.

Note that the negative of a (matroidal) M-convex set is also a
(matroidal) M-convex set, and the translation of a (matroidal)
M-convex set by an integral vector is also a (matroidal) M-convex set.
Therefore $d_{G}- \odotZ{B_{\rm out}\sp{\bullet}}$ is a matroidal M-convex set.
Clearly, a vector $m_{\rm in}$ is the in-degree vector of an
orientation $D$ of $G$ 
precisely if $d_{G}-m_{\rm in}$ is the out-degree vector of $D$.

We are interested in finding an orientation whose in-degree vector is
dec-min and whose out-degree vector is dec-min.  
This is equivalent to finding a member 
$m_{\rm in}$ of $\odotZ{B_{\rm in}\sp{\bullet}}$ 
for which the vector $m_{\rm out}:= d_{G}-m_{\rm in}$ 
is in the matroidal M-convex set $\odotZ{B_{\rm out}\sp{\bullet}}$.  
But this latter is equivalent to requiring that $m_{\rm in}$ 
is in the M-convex set $d_{G}- \odotZ{B_{\rm out}\sp{\bullet}}$.  
That is, the problem is equivalent to finding an element
of the intersection of the matroidal M-convex sets 
$\odotZ{B_{\rm in}\sp{\bullet}}$ and $d_{G}- \odotZ{B_{\rm out}\sp{\bullet}}$.  
This latter problem can
be solved by Edmonds matroid intersection algorithm \cite{Edmonds79}.

\section{In-degree constrained orientations of graphs}
\label{SCori2}

In this section we first describe an algorithm to find a dec-min
in-degree constrained orientation.  Second, we develop a complete
description of 
the set of dec-min in-degree constrained orientations,  
which gives rise to an algorithm to compute a
cheapest dec-min in-degree constrained orientation.

\subsection{Computing a dec-min in-degree constrained orientation}
\label{SCdecminfg}

Let $f:V\rightarrow {\bf Z}\cup \{-\infty \}$ be a lower bound function 
and $g:V\rightarrow {\bf Z}\cup \{+\infty \}$ an upper bound function 
for which $f\leq g$.  
We are interested in in-degree constrained orientations $D$ of $G$, 
by which we mean that $f(v)\leq \varrho_{D}(v)\leq g(v)$ for every $v\in V$.  
Such an orientation is called {\bf $(f,g)$-bounded}, 
and we assume that $G$ has such an orientation. 
By a well-known theorem (see, Theorem 2.3.5 in \cite{Frank-book}),
such an orientation exists 
if and only if $i_{G}\leq \widetilde g$ and $\widetilde f\leq e_{G}$.

As before, let $\odotZ{B_{G}}$ denote the M-convex set 
of the in-degree vectors of orientations of $G$, 
and let $\odotZ{B\sq_{G}}$ denote the intersection of $\odotZ{B_{G}}$ 
with the integral box $T(f,g)$.  
That is, $\odotZ{B\sq_{G}}$ is the set 
of in-degree vectors of $(f,g)$-bounded orientations of $G$.  
Let $D$ be an $(f,g)$-bounded orientation of $G$
with in-degree vector $m$.  
We denote the smallest tight set containing a node $t$ by 
$T\sq_{D}(t)$ ($=T\sq_m(t))$.  
By applying
Claim \ref{CLTm(u)}
to $\odotZ{B\sq_{G}}$, we obtain that
\begin{equation} \label{(Tsqdt)} 
 T\sq_{D}(t)
 = \begin{cases}
    \{ t  \} & \ \ \hbox{if}\ \ \ \varrho_{D}(t)=f(t), 
   \\
    T_{D}(t) - \{s:  \varrho_{D}(s) = g(s)\} & \ \ \hbox{if}\ \ \ \varrho_{D}(t)>f(t), 
 \end{cases} 
\end{equation}
implying that, in case $\varrho_{D}(t)>f(t)$, the set $T\sq_{D}(t)$ consists of
those nodes $s$ from which $t$ is reachable and 
for which $\varrho_{D}(s)<g(s)$.

Formula \eqref{(Tsqdt)} implies for distinct nodes $s$ and $t$ that the
vector $m':= m + \chi_{s} - \chi_{t}$ belongs to $\odotZ{B\sq_{G}}$ precisely
if there is an $st$-dipath 
(i.e. a dipath from $s$ to $t$) for which 
$\varrho_{D}(s)<g(s)$ and $\varrho_{D}(t)>f(t)$.  
We call such a dipath $P$
of $D$ {\bf reversible}.  Note that the dipath $P'$ of $D'$ obtained
by reorienting $P$ is reversible in $D'$.

If $P$ is a reversible $st$-dipath of $D$ for which 
$\varrho_{D}(t)\geq \varrho_{D}(s)+2$, 
then the orientation $D'$ is decreasingly smaller than $D$.  
We call such a dipath {\bf improving}.  
Therefore, reorienting an improving $st$-dipath 
corresponds to a 1-tightening step.  
Hence 
Theorem~\ref{THequidotone} implies the following extension of Theorem~\ref{pathrev}.

\begin{theorem} \label{Borfg} 
An $(f,g)$-bounded orientation $D$ of $G$ is
dec-min if and only if there is no improving dipath, that is, a dipath
from a node $s$ to a node $t$ for which 
$\varrho_{D}(t) \geq \varrho_{D}(s)+2$, $\varrho_{D}(s)<g(s)$, 
and $\varrho_{D}(t)>f(t)$.  
\finbox 
\end{theorem}

In Section \ref{SCbasicalg} we have presented an algorithm that computes
a dec-min element of an arbitrary M-convex set. 
By specializing it to $\odotZ{B\sq_{G}}$, 
we conclude that 
in order to construct a dec-min $(f,g)$-bounded orientation of $G$, 
one can start with an arbitrary $(f,g)$-bounded orientation, and then
reorient (currently) improving dipaths one by one, 
as long as such a dipath exists.  
As we pointed out after Theorem~\ref{pathrev}, 
after at most $\vert E\vert \sp{2}$ improving dipath reorientations, 
the algorithm terminates with a dec-min $(f,g)$-bounded orientation of $G$.

\medskip

\paragraph{Canonical chain and essential value-sequence 
for $(f,g)$-bounded orientations} \ 
In Section \ref{SCcompcanchain}, 
we indicated that Algorithm \ref{algo1} can immediately be applied to 
compute the canonical chain, 
the canonical partition, and the essential value-sequence 
belonging to the intersection $\odotZ{B\sq}$ of an arbitrary M-convex set $\odotZ{B}$ 
with an integral box $T(f,g)$.  

This algorithm needs only the original subroutine to compute $T_m(u)$
since, by Claim \ref{CLTm(u)}, 
$T\sq_m(u)$ is easily computable from $T_m(u)$.  
As we indicated above, in the special case of orientations, 
the corresponding sets $T_{D}(t)$ and $T\sq_{D}(t)$ are
immediately computable from $D$.  
Therefore this extended algorithm
can be used in the special case when we are interested in dec-min
$(f,g)$-bounded orientations of $G=(V,E)$.  
The algorithm starts with a dec-min $(f,g)$-bounded orientation $D$ of $G$
and outputs the canonical chain ${\cal C}\sq = \{C\sq_{1},\dots ,C\sq_{q}\}$, 
the canonical partition ${\cal P}\sq = \{S\sq_{1},\dots ,S\sq_{q}\}$, and 
the essential value-sequence $\beta \sq_{1}>\cdots >\beta \sq_{q}$.
In view of 
Theorem~\ref{THchardecmin},
we also define bounding functions $f\sp{*}$ and $g\sp{*}$ as
\begin{align*}
f\sp{*}(v) & := \beta \sq_{i} -1 \ \hbox{if}\ \ v\in S_{i} 
\qquad (i=1,\dots,q), 
\\
g\sp{*}(v) & := \beta \sq_{i} \phantom{{} -1}  \ \hbox{if}\ \ v\in S_{i} 
\qquad  (i=1,\dots ,q). 
\end{align*}
\noindent
We say that the small box 
\begin{equation} \label{(fgbounded)} 
T\sp{*}:=T(f\sp{*},g\sp{*}) 
\end{equation}
belongs to $\odotZ{B\sq_{G}}$.  
Clearly, $f\leq f\sp{*}$ and $g\sp{*} \leq g$, 
and hence $T(f\sp{*},g\sp{*})\subseteq T(f,g)$.
In Section \ref{SCoricheapest} 
below we assume that these data are available.

\begin{remark} \rm \label{Tspec} 
A special case of in-degree constrained
orientations is when we have a prescribed subset $T$ of $V$ and a
non-negative function $m_T:T\rightarrow {\bf Z}_+$ serving as an
in-degree specification on $T$, and we are interested in orientations
of $G$ for which $\varrho(v)=m_T(v)$ holds for every $t\in T$.  We
call such an orientation {\bf $T$-specified}. 
This notion will have applications in Section \ref{SCsemimatching}.
\finbox
\end{remark}

\subsection{Cheapest dec-min in-degree constrained orientations}
\label{SCoricheapest}

We are given a cost-function $c$ on the possible orientations of the edges of $G$ 
and our goal is to find a cheapest dec-min $(f,g)$-bounded orientation of $G$.  
This will be done with the help of a purely graphical description 
of the set of all dec-min $(f,g)$-bounded orientations, 
which is given in 
Theorem~\ref{THfgbounded}.

As a preparation, we derive the following claim as an immediate
consequence of 
the structural result stated in 
Theorem~\ref{THmatroid-eltolt}. 
Let $m$ be a dec-min element of an M-convex 
set $\odotZ{B}$ on ground-set $S$.
Suppose that
$m':=m + \chi_{s} - \chi_{t}$ is in $\odotZ{B}$ 
(that is, $s\in T_m(t)$).
Since $m$ is dec-min, $m(t)\leq m(s)+1$.  
If $m(t)=m(s)+1$, then $m'$ and $m$ are value-equivalent and 
hence $m'$ is also a dec-min element of $\odotZ{B}$.  
We say that $m'$ is obtained from $m$ by an {\bf elementary step}.

\begin{claim} \label{smallstep} 
Any dec-min element of $\odotZ{B}$ can be obtained from 
a given dec-min element $m$ by a sequence of at most $\vert S\vert $ elementary steps.  
\end{claim}

\Proof 
By Theorem~\ref{THmatroid-eltolt}, 
the set of dec-min elements of
$\odotZ{B}$ is a matroidal M-convex set in the sense that it
can be obtained from a matroid $M\sp{*}$ 
by translating the incidence vectors of the bases of $M\sp{*}$ 
by the same integral vector $\Delta\sp{*}$.  
A simple property of matroids is that any basis can be obtained from a given basis 
through a sequence of at most $\vert S\vert $ bases 
such that each member of the series can be obtained from 
the preceding one by taking out one element and adding a new one.
The corresponding change in the translated vector is exactly an elementary step.
\finbox 
\medskip

\begin{theorem} \label{THfgbounded} 
Let $G=(V,E)$ be an undirected graph admitting an $(f,g)$-bounded orientation.  
Let $\odotZ{B\sq_{G}}$ denote the M-convex set consisting of the in-degree vectors 
of $(f,g)$-bounded orientations of $G$, 
and let $T\sp{*}$ be the small box, 
belonging to $\odotZ{B\sq_{G}}$, 
as defined in \eqref{(fgbounded)}.
There is a chain $\cal Z$ of subsets of $V$ such that an
$(f,g)$-bounded orientation $D$ of $G$ is a dec-min $(f,g)$-bounded orientation 
if and only if 
$D$ is an orientation of $G$ whose in-degree vector belongs to 
$T\sp{*}$ and $\delta_D(Z)=0$ holds for each $Z\in {\cal Z}$. 
\end{theorem}

\Proof Let $D$ be a dec-min $(f,g)$-bounded orientation of $G$, and
let $m$ denote its in-degree vector.  Consider the canonical chain
${\cal C}\sq = \{C\sq_{1},\dots ,C\sq_{q}\}$, the canonical partition
${\cal P}\sq = \{S\sq_{1},\dots ,S\sq_{q}\}$, and the essential
value-sequence $\beta \sq_{1}>\cdots >\beta \sq_{q}$ belonging to
$\odotZ{B\sq_{G}}$.

For $i\in \{1,\dots ,q\}$, define $$F_{i}:=\{v:  v\in S\sq_{i}, f(v) =
\beta \sq_{i}\}.$$ Since $f(v)\leq m(v)\leq \beta \sq_{i}$ holds for every
element $v$ of $S\sq_{i}$, we obtain that $f(v)=m(v)=\beta \sq_{i}$ 
for $v \in F_{i}$.  
Note that $F_{i}$ does not depend on $D$.

\begin{claim} \label{nostpath} 
For every $h=1,\dots ,i$, there is no dipath $P$ 
from a node 
$s\in V-C\sq_{h}$
with $m(s)<g(s)$ 
to a node $t\in S\sq_{h}$ with $\beta \sq_{h}>f(t)$.  
\end{claim}

\Proof Suppose indirectly that there is such a dipath $P$.  
If $m(t) = \beta \sq_{h}$, then $P$ would be an improving dipath which is
impossible since $D$ is dec-min $(f,g)$-bounded.
Therefore $m(t)=\beta \sq_{h}-1$.  
But a property of the canonical partition is
that there is an element $t'$ of $S\sq_{h}-F_{h}$
for which $m(t')=\beta\sq_{h}$ 
and $t\in T\sq_{D}(t')$.  
This means that $t'$ is reachable from
$t$ in $D$, and therefore there is a dipath from $s$ to $t'$ in $D$
which is improving, a contradiction again.  
\finbox

\medskip

We are going to define a chain ${\cal Z}$ of subsets 
$Z_{1} \supseteq Z_{2} \supseteq \cdots \supseteq Z_{q} \ (=\emptyset )$ 
of $V$ with the help
of $D$, and will show that 
this chain actually does not depend on $D$.
Let 
\begin{equation} \label{(Zidef)} 
 Z_{i}:= \{ t: \, 
  \mbox{$t$ is reachable in $D$ from a node \, $s\in V-C\sq_{i}$
  \, with $\varrho_{D}(s)<g(s)$}   \}. 
\end{equation}
\noindent
Note that $Z_{i-1}\supseteq Z_{i}$ follows from the definition, where 
equality holds precisely if $\varrho_{D}(s) =g(s)$ for each $s\in S_{i}$.

\begin{lemma} \label{Zisame} 
Every dec-min $(f,g)$-bounded orientation defines the same family ${\cal Z}$.  
\end{lemma}

\Proof 
By Claim \ref{smallstep}, it suffices to prove that a single elementary step 
does not change ${\cal Z}$.  
An elementary step in $\odotZ{B\sq_{G}}$ corresponds 
to the reorientation of an $st$-dipath $P$ in $D$ 
where $s,t\in S\sq_{h}-F_{h}$, $m(t)=\beta \sq_{h}$ and 
$m(s)=\beta\sq_{h}-1$ hold 
for some $h\in \{1,\dots ,q\}$.  
We will show for $i\in \{1,\dots ,q\}$ 
that the reorientation of $P$ does not change $Z_{i}$.

If $h\leq i$, then Claim \ref{nostpath} implies that 
$Z_{i}\cap S\sq_{h} \subseteq F_{h}$.  
Since $\delta_{D}(Z_{i})=0$, the dipath $P$ is disjoint from $Z_{i}$, 
implying that reorienting $P$ does not affect $Z_{i}$.

Suppose now that $h\geq i+1$.  
Since reorienting $P$ results in a dec-min $(f,g)$-bounded orientation $D'$, 
we get that $m(s)+1\leq g(s)$ and hence $s\in Z_{i}-C\sq_{i}$.
Since $\delta_{D}(Z_{i})=0$, we obtain that $t\in Z_{i}-C\sq_{i}$.  
Since 
$\varrho_{D'}(t) = \varrho_{D}(t)-1 < g(t)$
and the set of nodes reachable from $s$ in $D$ is equal to the set of
nodes reachable from $t$ in $D'$, 
it follows that the reorientation of $P$ does not change $Z_{i}$.  
\finbox 
\medskip
\medskip

\begin{lemma} \label{LME0A0}
The chain $\cal Z$ defined above meets the requirements of the theorem. 
\end{lemma}

\Proof 
Consider first an arbitrary dec-min $(f,g)$-bounded orientation $D$ of $G$.  
Lemma \ref{Zisame} and the definition of $Z_i$ in
\eqref{(Zidef)} imply that no arc of $D$ can leave any $Z_i\in {\cal Z}$, that is, 
$\delta _D(Z_i)=0$.  
Furthermore, by applying Theorem \ref{THchardecmin}
to the present base-polyhedron $\odotZ {B\sq}$, 
we obtain that the in-degree vector of $D$ belongs to $T\sp{*}$.

Conversely, let $D$ be an orientation of $G$ for which 
$\delta_D(Z)=0$ holds for every $Z\in {\cal Z}$ and 
the in-degree vector of $D$ belongs to $T\sp{*}$, that is, 
\begin{equation*}
 \hbox{ $f\sp{*}(v) \leq \varrho_{D}(v) \leq g\sp{*}(v)$ for every $v\in V$.  }\
\end{equation*}
\noindent
Then $D$ is clearly $(f,g)$-bounded.
The following claim makes complete the proof of the lemma.

\begin{claim} \label{CLnoimproving}
There is no improving dipath in $D$.  
\end{claim}

\Proof 
Suppose, indirectly, that $P$ is an improving $st$-dipath, that is, 
a dipath from $s$ to $t$ such that 
$\varrho_{D}(t) \geq \varrho_{D}(s)+2$, 
$\varrho_{D}(t)>f(t)$, and $\varrho_{D}(s)<g(s)$.  
Suppose that $t$ is in $S\sq_{i}$ for some $i\in \{1,\dots ,q\}$.  
If $s$ is in $S_k\sq$ for some $k\in \{1,\dots ,q\}$, then
\[
 \beta \sq_k-1 \leq \varrho_{D}(s) \leq \varrho_{D}(t) -2\leq \beta
\sq_{i}-2,
\]
that is, $\beta \sq_k < \beta \sq_{i}$, and hence $k>i$,
implying that $s$ is in $V-C\sq_{i}$.  This and $\varrho_{D}(s)<g(s)$
imply that $s$ is in $Z_{i}$.  
Furthermore, $\beta \sq_{i}\geq \varrho_{D}(t)>f(t)$ implies that $t$ is not in $F_{i}$, 
and since 
$S_{i}\sq \cap Z_{i}\subseteq F_{i}$, 
we obtain that $t$ is not in $Z_{i}$.  
On the other hand, 
we must have $t \in Z_{i}$,
since there is a dipath 
from $s \in V-C\sq_{i}$ to $t$ and $\varrho_{D}(s)<g(s)$.
This is a contradiction.
\finbox
\medskip

By proving Claim \ref{CLnoimproving}, we have shown Lemma \ref{LME0A0}.
Thus the proof of 
Theorem~\ref{THfgbounded} is completed.
\finbox \finboxHere 
\medskip

\paragraph{Algorithm for computing a cheapest dec-min $(f,g)$-bounded orientation} \ 
First we compute a dec-min $(f,g)$-bounded orientation $D$ 
of $G$ with the help of the algorithm outlined in Section \ref{SCdecminfg}.  
Second, by applying the algorithm described in the same section, 
we compute the canonical chain and partition belonging to $\odotZ{B\sq_{G}}$ 
along with the essential value-sequence.  
Once these data are available, the sets $Z_{i}$
($i=1,\dots ,q$) defined in \eqref{(Zidef)} are easily computable.
Lemma \ref{Zisame} ensures that these sets $Z_{i}$ do not depend on the
starting dec-min $(f,g)$-bounded orientation $D$.  
Let $E_{0}$ be the
union of the set of edges connecting some $Z_{i}$ with $V-Z_{i}$, and
define the orientation $A_{0}$ of $E_{0}$ 
by orienting each edge between $Z_{i}$ and $V-Z_{i}$ toward $Z_{i}$.

Theorem~\ref{THfgbounded} implies that, once $E_{0}$ and its orientation
$A_{0}$ are available, the problem of computing a cheapest dec-min
$(f,g)$-bounded orientation of $G$ reduces to finding cheapest
in-degree constrained 
(namely, $(f\sp{*},g\sp{*})$-bounded) 
orientation of a mixed graph.  
We indicated already in Section \ref{cheapdecmin} 
that such a problem is easily 
solvable 
by the strongly polynomial
min-cost flow algorithm of Ford and Fulkerson in a digraph with identically 1 capacities.

\begin{remark} \rm \label{RMorivar}
In Section \ref{capori} we have considered the capacitated dec-min orientation problem
in the basic case where no in-degree constraints are imposed.
With the technique presented there, 
we can cope with the capacitated, min-cost, in-degree constrained variants as well.
Furthermore, 
the algorithms above can easily be extended, with a slight modification, 
to the case when one is interested in orientations of mixed graphs. 
\finbox
\end{remark}

\subsection{Dec-min $(f,g)$-bounded orientations minimizing the in-degree of $T$} 
\label{SCminm(T)}

In this section, we discuss a rather specific orientation problem in full detail.  
The reason is that this framework is an indispensable tool 
for solving in Section \ref{SCsemimatching} a common generalization
of several previously investigated resource allocation problems.  
Our approach permits us to manage algorithmically even the minimum cost
version of these problems.

One may consider $(f,g)$-bounded orientations of $G$ when the
additional requirement is imposed that the in-degree of a specified subset $T$ 
of nodes be as small as possible.  
We shall show that these orientations of $G$ can be described 
as $(f',g')$-bounded orientations
of a mixed graph arising from $G$ by orienting the edges 
between a certain subset $X_{T}$ 
of nodes and its complement $V-X_{T}$ toward $V-X_{T}$.

It is more comfortable, however, to show the analogous statement for a
general M-convex set $\odotZ{B'}(p)\subseteq {\bf Z}\sp V$ defined by a
(fully) supermodular function $p$ for which $\odotZ{B\sq}:=
\odotZ{B'}(p)\cap T(f,g)$ is non-empty.  (Here, instead of the usual
$S$, we use $V$ to denote the ground-set of the general M-convex set.
We are back at the special case of graph orientations when $p=i_G$.)
We assume that each of $p$, $f$, and $g$ is finite-valued.

Let $p\sq$ denote the unique (fully) supermodular function defining $B\sq$.  
This function can be expressed with the help of $p$, $f$, and $g$, as follows
(see, for example, Theorem 14.3.9 in \cite{Frank-book}):
\begin{equation}
p\sq(Y) = \max \{p(X) + \widetilde f(Y-X) - \widetilde g(X-Y):  \ X\subseteq V \} 
\qquad(Y \subseteq V).
\label{(pfgb)} 
\end{equation}
\noindent
As $B\sq$ is defined by the supermodular function $p\sq$
 (that is, $B\sq=B'(p\sq)$), we have
\begin{equation}
\min \{\widetilde m(T):  m\in \odotZ{B\sq} \}=p\sq(T) .
\label{(mminim)} 
\end{equation}
\noindent
This implies that the set of elements of $\odotZ{B\sq}$ minimizing $\widetilde m(T)$ 
is itself an M-convex set.  
Namely, it is the set of integral elements of the base-polyhedron arising from
$B\sq$ by taking its face defined by 
$\{m\in B\sq:  \widetilde m(T)=p\sq(T)\}$.  
The next theorem shows how this M-convex set can be described 
in terms of $f$, $g$, and $p$, without referring to $p\sq$.

\begin{theorem}  \label{minmT} 
There is a box $T(f',g')\subseteq T(f,g)$ and a subset 
$X_{T}\subseteq V$ such that an element $m\in \odotZ{B\sq}$
minimizes $\widetilde m(T)$ 
if and only if 
$\widetilde m(X_{T})=p(X_{T})$ and $m\in \odotZ{B}\cap T(f',g')$.  
\end{theorem}

\Proof Let $X_{T}$ be a set maximizing the right-hand side of \eqref{(pfgb)}
for the definition of $p\sq(Y)$.

\begin{claim}  \label{opcri} 
An element $m\in \odotZ{B\sq}$ is a minimizer of
the left-hand side of \eqref{(mminim)} 
if and only if 
the following three optimality criteria hold:
\begin{align*}
& \widetilde m(X_{T})=p(X_{T}),
\\
& v\in T-X_T \quad \mbox{\rm implies} \quad m(v)=f(v) ,
\\ 
& v\in X_T-T \quad \mbox{\rm implies} \quad m(v)=g(v).
\end{align*}
\end{claim}

\Proof 
For any $m\in \odotZ{B\sq}$ and $X\subseteq V$, we have
$\widetilde m(T) = \widetilde m(X) + \widetilde m(T-X) - \widetilde m(X-T) 
 \geq p(X) + \widetilde f(T-X) - \widetilde g(X-T)$.  
Here we have equality if and only if 
$\widetilde m(X) =p(X)$, $\widetilde m(T-X)= \widetilde f(T-X)$, 
and $\widetilde m(X-T)= \widetilde g(X-T)$, 
implying the claim.  
\finbox

\medskip 
Define $f'$ and $g'$ as follows:
\begin{align}
 f'(v) &:= 
 \begin{cases} 
   g(v) & \ \ \hbox{if}\ \ \ v\in X_{T}-T, \cr
   f(v) & \ \ \hbox{if}\ \ \ v\in V- (X_{T}-T), 
 \end{cases} 
\label{(f'def)}
\\
g'(v) &:= 
 \begin{cases} 
  f(v) & \ \ \hbox{if}\ \ \ v\in T-X_{T}, \cr
  g(v) & \ \ \hbox{if}\ \ \ v\in V- (T-X_{T}).  
 \end{cases}
\label{(g'def)} 
\end{align}
\noindent
The claim implies that $T(f',g')$ and $X_{T}$ meet the requirement of the theorem.  
\finbox \finboxHere \medskip

As the set of elements of $\odotZ{B\sq}$ minimizing $\widetilde m(T)$
is itself an M-convex set, all the algorithms developed earlier can be
applied once we are able to compute set $X_{T}$ occurring in Theorem~\ref{minmT}.  
(By definitions \eqref{(f'def)} and \eqref{(g'def)}, 
$X_{T}$ immediately determines $f'$ and $g'$).

The following straightforward algorithm computes an element 
$m\in \odotZ{B\sq}$ minimizing the left-hand side of \eqref{(mminim)}
and a subset $X_{T}$ maximizing the right-hand side of \eqref{(pfgb)}.  
Start with an arbitrary element $m\in \odotZ{B\sq}$.  
By an improving step we mean the change of 
$m$ to $m':=m + \chi_{s} - \chi_{t}$ for some elements
$s\in V-T, t\in T$ for which $m(s)<g(s)$, $m(t)>f(t)$,  and $s\in T_{m}(t)$, 
where $T_{m}(t)$ is the smallest $m$-tight set 
(with respect to $p$) containing $t$.  
Clearly, $m'\in \odotZ{B\sq}$, and $\widetilde m'(T)=\widetilde m(T)-1$.  
The algorithm applies improving steps as long as possible.  
When no more improving step exists, the set 
$X_{T}:= \cup (T_{m}(t):  t\in T, m(t)>f(t))$ 
meets the three optimality criteria.  
The algorithm is polynomial if $\vert p(X)\vert $ is
bounded by a polynomial of $\vert V\vert $.

By applying Theorem~\ref{minmT} to the special case of $p=i_G$, we obtain the following.

\begin{corollary}  \label{grafminmt} 
Let $G=(V,E)$ be a graph admitting an
$(f,g)$-bounded orientation.  
There is a box $T(f',g')\subseteq T(f,g)$ 
and a subset $X_{T}\subseteq V$ such that an $(f,g)$-bounded
orientation of $G$ minimizes the in-degree of $T$ if and only if $D$
is an $(f',g')$-bounded orientation for which $\varrho _D(X_{T})=0$.
\finbox
\end{corollary}

In this case, the algorithm above to compute $X_{T}$ starts with an
$(f,g)$-bounded orientation $D$ of $G$, whose in-degree vector is
denoted by $m$.  As long as there is an $st$-dipath $P$ with $s\in
V-T, t\in T, m(s)<g(s)$, and $m(t)>f(t)$, reorient $P$.  When no such
a dipath exists anymore, the set $X_{T}$ of nodes from which a node
$t\in T$ with $m(t)>f(t)$ is reachable in $D$, along with the bounding
functions $f'$ and $g'$ defined in \eqref{(f'def)} and in
\eqref{(g'def)}, meet the requirement in the corollary.  \medskip

\paragraph{Minimum cost version} \ 
It follows that, in order to
compute a minimum cost dec-min $(f,g)$-bounded orientation for which
the in-degree of $T$ is minimum, we can apply the algorithm described
in Section~\ref{SCoricheapest} for the mixed graph
obtained from $G$ by orienting each edge between $X_{T}$ and $V-X_{T}$
toward $V-X_{T}$.

\begin{remark} \rm 
Instead of a single subset $T$ of $V$, we may consider a
chain ${\cal T}$ of subsets $T_{1}\subset T_{2}\subset \cdots \subset T_{h}$ of $V$.  
Then $\cal T$ defines a face $B\sq_{\rm face}$ of the base-polyhedron $B\sq$.  
Namely, an element $m$ of $B\sq$ belongs to
$B\sq_{\rm face}$ precisely if 
$\widetilde m(T_{i})=p\sq(T_{i})$
for each $i\in \{1,\dots ,h\}$.  
This implies that the integral elements of
$B\sq_{\rm face}$ simultaneously minimize $\widetilde m(T_{i})$ 
for each $i\in \{1,\dots ,h\}$ (over the elements of $\odotZ{B\sq}$).
Therefore, we can consider $(f,g)$-bounded orientations of $G$ with
the additional requirement that each of the in-degrees of 
$T_{1},T_{2}, \dots ,T_{h}$ 
is (simultaneously) minimum.  Corollary \ref{grafminmt}
can be extended to this case, implying that we have an algorithm to
compute a minimum cost dec-min $(f,g)$-bounded orientation of $G$ that
simultaneously minimizes the in-degree of each member of the chain
$\{T_{1}, T_{2}, \dots , T_{h}$\}. 
\finbox
\end{remark}

\subsection{Application in resource allocation: semi-matchings}
\label{SCsemimatching}

For a general M-convex set $\odotZ{B}$, 
it is known (see Theorem \ref{THdecmin=squaresum} in Section \ref{SCbasicalg})
that an element $m$ of $\odotZ{B}$ is dec-min
if and only if $m$ is a minimizer of the square-sum of the components over $\odotZ{B}$.
Therefore the corresponding equivalences hold 
in the special case of in-degree constrained 
(in particular, $T$-specified)
 orientations of undirected graphs.

As an application of this equivalence, 
we show first how a result of Harvey et al.~\cite{HLLT} 
concerning a resource allocation problem follows immediately.  
They introduced the notion of a semi-matching of 
a simple bipartite graph $G=(S,T;E)$ as a subset $F$ of edges 
for which $d_{F}(t)=1$ holds for every node $t\in T$, 
and solved the problem of finding a semi-matching $F$ for which
$\sum [d_{F}(s) (d_{F}(s)+1):  s\in S]$ is minimum.
(Here $d_{F}(v)$ denotes the number of edges in $F$ for which $v$ is an end-node.)
The problem was motivated by practical applications 
in the area of resource allocation in computer science.  
Note that
\begin{align*}
& \sum [d_{F}(s) (d_{F}(s)+1):  s\in S] 
 = \sum [d_{F}(s)\sp{2} :s\in S] + \sum [d_{F}(s):  s \in S] 
\\ &
= \sum [d_{F}(s)\sp{2} :s\in S] + \vert F\vert 
= \sum [d_{F}(s)\sp{2}:s\in S] + \vert T\vert ,
\end{align*}
and therefore the problem of Harvey et al.~is equivalent to
finding a semi-matching $F$ of $G$ that minimizes the square-sum of
degrees in $S$.

By orienting each edge in $F$ toward $S$ and each edge in $E-F$ toward $T$, 
a semi-matching can be identified with the set of arcs directed toward $S$ 
in an orientation of $G=(S,T;E)$ in which the out-degree of every node $t\in T$ is 1 
(that is, $\varrho(t)=d_{G}(t)-1$), and
$d_{F}(s)=\varrho(s)$ for each $s\in S$.  
Since $\varrho(t)$ for $t\in T$ is the same in these orientations, 
it follows that the total sum of $\varrho(v)\sp{2}$ over $S\cup T$ 
is minimized precisely if 
$\sum [\varrho(s)\sp{2}: s\in S] = \sum [d_{F}(s)\sp{2}:s\in S]$ is minimized.
Therefore the semi-matching problem of Harvey et al. is nothing but a
special dec-min $T$-specified orientation problem.
Note that not only semi-matching problems can be managed
with graph orientations, but conversely, 
an orientation of a graph $G=(V,E)$ can also be interpreted 
as a semi-matching of the bipartite
graph obtained from $G$ by 
subdividing each edge by a new node (where subdividing an edge
$e=uv$ formally means that we replace $e$ by a path $(uz_e,z_ev)$ of
length two with a new node $z_e$).  
This implies, for example, that the algorithm of Harvey et al.~to 
compute a semi-matching minimizing 
$\sum [d_F(v)\sp{2}: v\in S]$ 
is able to compute an orientation of a graph $G$ 
for which $\sum [\varrho(v)\sp{2}:v\in S]$ is minimum.  
Furthermore, an orientation of a hypergraph means that 
we assign an element of each hyper-edge $Z$ to $Z$ as its head.  
In this sense, semi-matchings of bipartite graphs and
orientations of hypergraphs are exactly the same.  
Several graph orientation results have been extended to hypergraph orientation, 
for an overview, see, e.g.  \cite{Frank-book}.

Bokal et al.~\cite{BBJ} extended the results to subgraphs of $G$
meeting a more general degree-specification on $T$ when, rather than
the identically 1 function, one imposes an arbitrary
degree-specification $m_{T}$ on $T$ satisfying 
$0\leq m_{T}(t)\leq d_{G}(t)$ $(t\in T)$.  
The same orientation approach applies in this more general setting.  
We may call a subset $F$ of edges an {\bf
$m_{T}$-semi-matching} if $d_{F}(t)= m_{T}(t)$ for each $t\in T$.  
The extended resource allocation problem is to find an $m_{T}$-semi-matching
$F$ that minimizes $\sum [ d_{F}(s)\sp{2}:  s\in S]$.  
This is equivalent to finding a $T$-specified orientation of $G$ 
for which the square-sum
of the in-degrees is minimum and the in-degree specification 
in $t\in T$ is 
$m_{T}'(t):  = d_{G}(t) - m_{T}(t)$.  
Therefore this extended resource allocation problem is equivalent 
to finding a dec-min $T$-specified orientation of $G$.

The same orientation approach, when applied to in-degree constrained
orientations, allows us to extend the $m_{T}$-semi-matching problem when
we have upper and lower bounds imposed on the nodes in $S$.  
This may be a natural requirement in practical applications
where the elements of $S$ correspond to available resources (e.g. computers), 
the elements of $T$ correspond to users, and we are interested in a fair
($=$ dec-min $=$ square-sum minimizer) distribution 
(=$m_{T}$-semi-matchings) of the resources when the load (or burden) of
each resource is requested to meet a specified upper and/or lower bound. 
Note that in the resource allocation framework, the degree $d_{F}(s)$ of
node $s\in S$ may be interpreted as the burden of $s$, and hence a
difference-sum minimizer semi-matching minimizes the total sum of
burden-differences.

Katreni{\v c} and Semani{\v s}in \cite{Katrenic13}
investigated the problem of finding a dec-min 
``maximum $(f,g)$-semi-matching'' problem where there is a lower-bound
function $f_{T}$ on $T$ and an upper bound function $g_{S}$ on $S$ 
(in the present notation) and one is interested in maximum cardinality
subgraphs of $G$ meeting these bounds.  
They describe an algorithm to compute a dec-min subgraph of this type.  
With the help of the orientation model discussed in Section \ref{SCminm(T)}
(where, besides the in-degree bounds on the nodes, 
the in-degree of a specified subset $T$ was requested to be minimum), 
we have a strongly polynomial algorithm for an extension of the model of
\cite{Katrenic13} when there may be upper and lower bounds on both $S$ and $T$.  
Actually, even the minimum cost version of this problem 
was solved in Section \ref{SCminm(T)}.

In another variation, we also have degree bounds $(f_{S},g_{S})$ on $S$
and $(f_{T},g_{T})$ on $T$, 
but we impose an arbitrary positive integer $\gamma $ 
for the cardinality of $F$.  
We consider degree-constrained subgraphs $(S,T;F)$ of $G$ 
for which $\vert F\vert =\gamma $, and want
to find such a subgraph for which 
$\sum [d_F(s)\sp{2} :  s\in S]$
is minimum.  
(Notice the asymmetric role of $S$ and $T$.)  
This is equivalent to finding an in-degree constrained orientation $D$ of $G$
for which 
$\varrho_D(S)=\gamma $ and 
$\sum [\varrho_D(s)\sp{2}: s\in S]$ 
is minimum.  
Here the corresponding in-degree bound $(f,g)$
on $S$ is the given $(f_{S},g_{S})$ 
while $(f,g)$ on $T$ is defined for $t\in T$ by
\[
 f(t):= d_{G}(t) -g_{T}(t) \quad \mbox{and} \quad g(t):= d_{G}(t) -f_{T}(t).
\]

Let $B$ denote the base-polyhedron spanned by the in-degree vectors of
the degree-constrained orientations of $G$.  Then the restriction of
$B$ to $S$ is a g-polymatroid $Q$.  By intersecting $Q$ with the hyperplane 
$\{x:  \widetilde x(S)=\gamma \}$, 
we obtain an integral base-polyhedron $B_{S}$ in ${\bf R}\sp{S}$, 
and then the elements of $\odotZ{B_{S}}$ are exactly 
the in-degree vectors of the requested orientations restricted to $S$.  
That is, the elements of $\odotZ{B_{S}}$ are the restriction of 
the degree-vectors of the requested subgraphs of $G$ to $S$.  
Since $B_{S}$ is a base-polyhedron, a dec-min element of $\odotZ{B_{S}}$ 
will be a solution to our minimum degree-square sum problem.

\medskip

We briefly indicate that a capacitated version of the semi-matching
problem can also be formulated as 
a dec-min in-degree constrained and capacitated orientation problem 
(cf., Section~\ref{capori} and Remark~\ref{RMorivar}).
Let $G=(S,T;E)$ be again a bipartite graph, $\gamma $ a positive integer, 
and $f_{V}$ and $g_{V}$ integer-valued bounding functions on
$V:=S\cup T$ for which $f_{V}\leq g_{V}$.  
In addition, an integer-valued capacity function $g_E$ 
is also given on the edge-set $E$, and we are interested 
in finding a non-negative integral vector $z:E\rightarrow {\bf Z}_{+}$ 
for which 
$\widetilde z(E)=\gamma $, $z\leq g_E$ and
$f_{V}(v)\leq d_z(v)\leq g_{V}(v)$
for every $v\in V$.  
(Here $d_z(v):=\sum [z(uv):  uv\in E]$.) 
We call such a vector {\bf feasible}.  
The problem is to find a feasible vector $z$ whose degree vector 
restricted to $S$ (that is, the vector $(d_z(s):  s\in S)$ is decreasingly minimal.

By replacing each edge $e$ with $g_E(e)$ parallel edges, it follows
from the uncapacitated case above that the vectors 
$\{(d_z(s):  s\in S) :  \mbox{$z$ is a feasible integral vector} \}$ 
form an M-convex set.  
In this case, however, the basic algorithm is not necessarily polynomial
since the values of $g_E$ may be large.  
Therefore we need the general
strongly polynomial algorithm described in Section \ref{strong.pol}.
 In this case the general Subroutine
\eqref{(ND.routine)} can be realized via max-flow min-cut computations.

\paragraph{Minimum edge-cost dec-min semi-matchings} \ 
Harada et al.~\cite{HOSY07} developed an algorithm to solve the minimum
edge-cost version of the original semi-matching problem of Harvey et al.~\cite{HLLT}.  
As the dec-min in-degree bounded orientation problem covers all the
extensions of semi-matching problems mentioned above, the minimum 
edge-cost version of these extensions can also be solved with the strongly
polynomial algorithms developed in Section \ref{SCoricheapest}
for minimum cost dec-min in-degree bounded orientations.

\medskip

We close this section with some historical remarks.  
The problem of Harvey et al.~is closely related to earlier investigations in the
context of minimizing a separable convex function over 
(integral elements of) a base-polyhedron.  
For example, Federgruen and Groenevelt \cite{Fed-Gro86} provided 
a polynomial time algorithm in 1986.  
Hochbaum and Hong \cite{Hochbaum-Hong} in 1995 developed 
a strongly polynomial algorithm;
their proof, however, included a technical gap,
which was fixed by Moriguchi, Shioura, and Tsuchimura \cite{MST} in 2011.  
For an early book on resource allocation, see the
one by Ibaraki and Katoh \cite{Ibaraki-Katoh.88} 
while three more recent surveys are due to 
Katoh and Ibaraki \cite{Katoh-Ibaraki.98} from 1998, 
to Hochbaum \cite{Hoc07} from 2007, and 
to Katoh, Shioura, and Ibaraki \cite{KSI} from 2013.
Algorithmic aspects of minimum degree square-sum problems for general
graphs were discussed by Apollonio and Seb{\H o} \cite{Apollonio-Sebo}.

\section{Orientations of graphs with edge-connectivity requirements}
\label{SCori3}

In this section, we investigate various edge-connectivity requirements
for the orientations of $G$.  
The main motivation behind these
investigations is a conjecture of Borradaile et al.~\cite{BIMOWZ} 
on decreasingly minimal strongly connected orientations.  
Our goal is to prove their conjecture in a more general form.

\subsection{Strongly connected orientations}
\label{SCstrori}

Suppose that $G$ is 2-edge-connected, implying that it has a strong
orientation by a theorem of Robbins \cite{Robbins}.  
We are interested in dec-min strong orientations, meaning 
that the in-degree vector is decreasingly minimal over the strong orientations of $G$.  
This problem of Borradaile et al.~was motivated by a practical application
concerning optimal interval routing schemes.

Analogously to Theorem~\ref{pathrev}, they described a natural way to
improve a strong orientation $D$ to another one whose in-degree
vector is decreasingly smaller.  
Suppose that there are two nodes $s$ and $t$ for which 
$\varrho(t)\geq \varrho(s)+2$ 
and there are two edge-disjoint dipaths from $s$ to $t$ in $D$.  
Then reorienting an arbitrary $st$-dipath of $D$ results in another strongly connected
orientation of $D$ which is clearly decreasingly smaller than $D$.

Borradaile et al. \cite{BIMOWZ} 
conjectured the truth of the converse (and this
conjecture was the starting point of our investigations).  The next
theorem states that the conjecture is true.

\begin{theorem}   \label{proved.conj} 
A strongly connected orientation $D$ of $G=(V,E)$ is decreasingly minimal 
if and only if 
there are no two arc-disjoint $st$-dipaths in $D$ 
for nodes $s$ and $t$ with $\varrho(t)\geq \varrho(s)+2$.
\end{theorem}

\Proof 
Suppose first that there are nodes $s$ and $t$ with 
$\varrho(t)\geq \varrho(s)+2$ 
such that there are two arc-disjoint $st$-dipaths of $D$.  
Let $P$ be any $st$-dipath in $D$ and let $D'$
denote the digraph arising from $D$ by reorienting $P$.  
Then $D'$ is strongly connected, 
since if it had a node-set $Z$ ($\emptyset \subset Z\subset V$) 
with no entering arcs, then $Z$ must be a $t\overline{s}$-set
and $P$ enters $Z$ exactly once.  
But then 
$0 = \varrho_{D'}(Z) =\varrho_{D}(Z)-1 \geq 2-1=1$, 
a contradiction.  
Therefore $D'$ is indeed strongly connected
and its in-degree vector is decreasingly 
smaller than that of $D$.

To see the non-trivial part, define a set-function $p_{1}$ as follows:
\begin{equation*} 
p_{1}(X):= 
\begin{cases}
 0 & \ \ \hbox{if}\ \ \ X=\emptyset,
\cr
\vert E\vert & \ \ \hbox{if}\ \ \ X=V ,
\cr 
i_{G}(X)+1 & \ \ \hbox{if}\ \ \ \emptyset \subset X\subset V. 
\end{cases} 
\end{equation*}
Then $p_{1}$ is crossing supermodular and hence $B_{1}:=B'(p_{1})$ is a base-polyhedron.

\begin{claim}   \label{mstrongori} 
An integral vector $m$ is the in-degree vector of a strong orientation of $G$
if and only if 
$m$ is in $\odotZ{B_{1}}$.
\end{claim}  

\Proof 
If $m$ is the in-degree vector of a strong orientation of $G$,
then $\widetilde m(V)= \vert E\vert = p_{1}(V)$, \ 
$\widetilde m(\emptyset) = 0 = p_{1}(\emptyset )$, and
\[
\widetilde m(Z) = \sum [\varrho(v):v\in Z] 
= \varrho(Z) + i_{G}(Z) \geq 1 + i_{G}(Z)= p_{1}(Z)
\]
for $\emptyset \subset Z\subset V$, that is, $m\in \odotZ{B_{1}}$.

Conversely, let $m\in \odotZ{B_{1}}$.  
Then $m\in B_{G}$ and hence by Claim \ref{oribase}, 
$G$ has an orientation $D$ with in-degree vector $m$.  
We claim that $D$ is strongly connected.  
Indeed,
\[
\varrho(Z)
 = \sum [\varrho(v):v\in Z] - i_{G}(Z) 
 = \widetilde m(Z) - i_{G}(Z) 
 \geq p_{1}(Z)-i_{G}(Z) =1
\]
whenever $\emptyset \subset Z\subset V$.  
\finbox

\begin{claim} \label{strongreori} 
Let $D$ be a strong orientation of $G$ with in-degree vector $m$.  
Let $t$ and $s$ be nodes of $G$.  
The vector $m':=m + \chi_{s} - \chi_{t}$ is in $B_{1}$ 
if and only if 
$D$ admits two arc-disjoint dipaths from $s$ to $t$.  
\end{claim}

\Proof 
$m'\in B_{1}$ holds precisely if there is no $t\overline{s}$-set $X$
which is  $m$-tight
 with respect to $p_{1}$, 
that is, 
$\widetilde m(X) = i_{G}(X)+1$.  
Since $\varrho(Y) + i_{G}(Y) = \sum [\varrho(v):v\in Y]= \widetilde m(Y)$ 
holds for any set $Y\subset V$, the tightness of $X$
(that is, $\widetilde m(X)= i_{G}(X)+1$) 
is equivalent to requiring that $\varrho(X)=1$.  
Therefore $m'\in B_{1}$ 
if and only if 
$\varrho(Y) > 1$
holds for every $t\overline{s}$-set $Y$, 
which is, by Menger's theorem, equivalent to 
the existence of two arc-disjoint $st$-dipaths of $D$.
\finbox 
\medskip

By Theorem~\ref{THequidotone},
$m$ is a dec-min element of $\odotZ{B_{1}}$ 
if and only if 
there is no 1-tightening step for $m$.  
By Claim \ref{strongreori} 
this is just equivalent to the condition in the theorem 
that there are no two arc-disjoint $st$-dipaths in $D$ 
for nodes $s$ and $t$ for which $\varrho(t)\geq \varrho(s)+2$.  
\finbox \finboxHere

\medskip

We remark that a recent work of Zhou and Hou \cite{Zhou+Hou}
describes a proof of the conjecture of Borradaile et al.~that 
uses purely graph-theoretical concepts and does not rely on any knowledge
of base-polyhedra or supermodular functions.  
An immediate consequence of Claim \ref{mstrongori} and 
Theorem~\ref{THequidotone} is the following.

\begin{corollary} 
A strong orientation of $G$ is dec-min 
if and only 
if it is inc-max.  
\finbox 
\end{corollary}

We indicated in Section \ref{SCdecminfg}
how in-degree constrained dec-min orientations can
be managed due to the fact that the intersection of an integral
base-polyhedron $B$ with an integral box $T$ is an integral base-polyhedron.  
The same approach works for degree-constrained strong orientations.  
For example, in this case dec-min and inc-max
again coincide and one can formulate the in-degree constrained
version of Theorem~\ref{proved.conj}.
In Section~\ref{SChighori}, 
we overview more general cases.

\subsection{Counterexample for mixed graphs}
\label{SCstroricntex}

Although Robbins' theorem on strong orientability of undirected graphs
easily extends to mixed graphs, 
as was pointed out by Boesch and Tindell \cite{BT}, 
it is not true anymore that a decreasingly minimal strong orientation 
of a mixed graph is always increasingly maximal.
Actually, one may consider two natural variants.

\begin{figure}
\centering
\includegraphics[width=0.8\textwidth,clip]{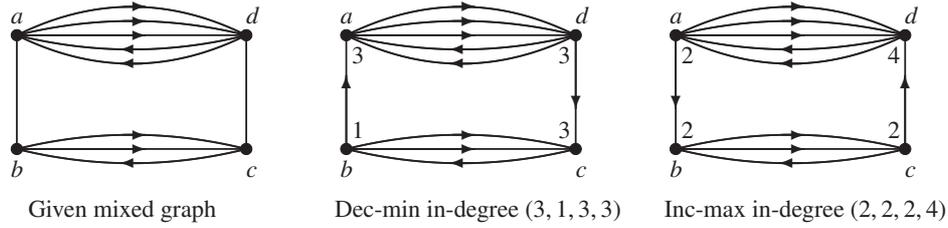}
\caption{Dec-min and inc-max strong orientations of a mixed graph may differ}
\label{FGMixG4}
\end{figure}

In the first one, decreasing minimality and increasing maximality
are defined for the total in-degree vector of the directed graph obtained from
the initial mixed graph orienting its undirected edges.  
Consider the mixed graph $M=(V,A+E)$ 
in the left of Figure \ref{FGMixG4}.
Here $V=\{a,b,c,d\}$ while the set $E$ of undirected edges 
(to be oriented) has just two elements:  $ab$ and $dc$.  
There are two strong orientations of $M$.  
In the first one (the middle digraph in Figure \ref{FGMixG4}), 
these are the arcs $ba$ and $dc$, 
in which case the total in-degree vector is $(3,1,3,3)$.  
In the second one (the digraph in the right 
of Figure \ref{FGMixG4}), the orientations of the elements of
$E$ are $ab$ and $cd$, in which case the total in-degree vector is $(2,2,2,4)$.  
Now $(3,1,3,3)$ is dec-min while $(2,2,2,4)$ is inc-max.

\begin{figure}
\centering
\includegraphics[width=0.35\textwidth,clip]{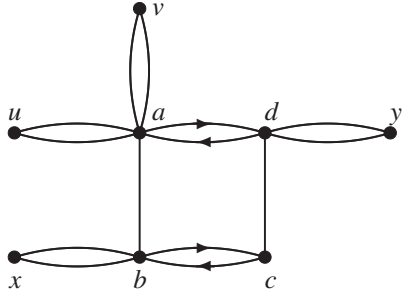}
\caption{A mixed graph for the second variant} 
\label{FGMixG8}
\end{figure}

In the second variant, we are interested in the in-degree vector of
the digraph obtained by orienting the originally undirected part $E$.
For this version, the counterexample mixed graph is demonstrated in Figure \ref{FGMixG8}. 
In any strong orientations of $M$, parallel edges in $E$ must be oriented oppositely.  
Therefore the contribution of the actual orientations of these 
four pairs of parallel undirected edges 
(in the order of $a, b, c, d, u, v, x, y$) 
is $(2, 1, 0, 1, 1, 1, 1, 1)$.

Therefore there are essentially two distinct strong orientations of $M$.  
In the first one, the undirected edges $ab, cd$ are oriented as $ba, dc$, 
while in the second one the undirected edges $ab, cd$ are oriented as $ab, cd$.  
Hence the in-degree vector of the first strong
orientation corresponding to the orientation of $G$ 
(in the order of $a, b, c, d, u, v, x, y$)
 is \ $(3, 1, 1, 1, 1, 1, 1, 1)$.  
The in-degree vector of second strong orientation corresponding to the
orientation of $G$ is $(2,2,0,2,1,1,1,1)$.  
The first vector is inc-max while the second vector is dec-min.

\medskip

These examples give rise to the question:  
what is behind the phenomenon that 
while dec-min and inc-max coincide for strong orientations of undirected graphs, 
they differ for strong orientations of mixed graph?  
The explanation is, as we pointed out earlier, that
for an M-convex set the two notions coincide and the set of 
in-degree vectors of strong orientations of an 
undirected graph is an M-convex set, 
while the corresponding set for a mixed graph is, in general, 
not an M-convex set.  
It is actually the intersection of two M-convex sets.  
An algorithm for computing a dec-min element of the intersection of two M-convex sets
is described in \cite{FM20partD}.

\medskip

\subsection{Higher edge-connectivity}
\label{SChighori}

An analogous approach works in a much more general setting.  
We say that a digraph 
{\bf covers}
 a set-function $h$ if $\varrho(X)\geq h(X)$
holds for every set $X\subseteq V$.  
The following result was proved in \cite{FrankJ4}.

\begin{theorem}[\cite{FrankJ4}]  \label{hori} 
Let $h$ be a finite-valued, non-negative
crossing supermodular function with $h(\emptyset )=h(V)=0$.  A graph
$G=(V,E)$ has an orientation covering $h$ 
if and only if
\[
 e_{\cal P} \geq \sum_{i=1}\sp q h(V_{i}) \quad \hbox{\rm and}\quad 
 e_{\cal P} \geq \sum_{i=1}\sp q h(V-V_{i}) 
\]
hold for every partition ${\cal P}=\{V_{1},\dots ,V_{q}\}$ of $V$, 
where $e_{\cal P}$ denotes the number of edges connecting
distinct parts of $\cal P$.  
\finbox 
\end{theorem}

This theorem easily implies the classic orientation result of
Nash-Williams \cite{NWir} stating that a graph $G$ has a
$k$-edge-connected orientation precisely if $G$ is $2k$-edge-connected.  
Namely, by defining $h$ to be identically equal to $k$ on non-empty proper
subsets of $V$ while $h(\emptyset ):=h(V):=0$, and observing that
$e_{\cal P}= \sum [d_G(X)/2:X\in {\cal P}]$, we can see that both
inequalities in the theorem are equivalent to 
$\sum [d_G(X)/2:X\in {\cal P}]\geq kq$.  
But this latter inequality follows from the
$(2k)$-edge-connectivity of $G$ since 
$\sum [d_G(X)/2:X\in {\cal P}]\geq \sum [k:  X\in {\cal P}] = kq$.

Even more, call a digraph $(k,\ell)$-edge-connected ($\ell\leq k$) 
(with respect to a root-node $r_{0}$) 
if $\varrho(X)\geq k$ whenever 
$\emptyset \subset X\subseteq V-r_{0}$ and 
$\varrho(X)\geq \ell$ whenever $r_{0}\in X\subset V$.  
(By Menger's theorem, 
$(k,\ell)$-edge-connectedness is equivalent to requiring that 
 there are $k$ arc-disjoint dipaths from $r_{0}$ to every node and
there are $\ell$ arc-disjoint dipaths from every node to $r_{0}$.)
Then Theorem~\ref{hori} implies:

\begin{theorem}  \label{kellori} 
A graph $G=(V,E)$ has a $(k,\ell)$-edge-connected orientation 
if and only if
\[
 e_{\cal P} \geq k(q-1) + \ell 
\]
holds for every $q$-partite partition $\cal P$ of $V$.  
\finbox 
\end{theorem}  

Note that an even more general special case of Theorem~\ref{hori} 
can be formulated to characterize graphs
admitting in-degree constrained and $(k,\ell)$-edge-connected orientations.

It is important to emphasize that however general Theorem~\ref{hori} is, 
it does not say anything about strong orientations of mixed graphs.  
In particular, it does not imply the pretty but easily provable theorem 
of Boesch and Tindell \cite{BT}.  
The problem of finding decreasingly minimal in-degree constrained 
$k$-edge-connected orientation of mixed graphs 
can be solved as a special case of decreasing minimization 
over the intersection of two M-convex sets.

The next lemma shows why the set of in-degree vectors of orientations
of $G$ covering the set-function $h$ appearing in Theorem~\ref{hori}
is an M-convex set, ensuring in this way the possibility 
of applying
the results on decreasing minimization over M-convex sets to general
graph orientation problems.

\begin{lemma}  \label{hbase} 
An orientation $D$ of $G$ covers $h$ 
if and only if 
its in-degree vector $m$ is in the base-polyhedron $B=B'(p)$,
where $p:= h+i_{G}$ is a crossing supermodular function.  
\end{lemma}

\Proof 
Suppose first that $m$ is the in-degree vector of a digraph covering $h$.  
Then $h(X)\leq \varrho(X)= \widetilde m(X) - i_{G}(X)$
for $X\subset V$ and 
$h(V) = 0 = \varrho(V)= \widetilde m(V) - i_{G}(V)$, 
that is, $m$ is indeed in $B$.

Conversely, suppose that $m\in B$.  
Since $h$ is finite-valued and non-negative, we have 
$\widetilde m(X) \geq p(X)\geq i_{G}(X)$ 
for
$X\subset V$ and $\widetilde m(V)=i_{G}(V)$ 
and hence, by the Orientation lemma, 
there is an orientation $D$ of $G$ with in-degree vector $m$.  
Moreover, this digraph $D$ covers $h$ since 
$\varrho_{D}(X) = \widetilde m(X) - i_{G}(X) \geq p(X) - i_{G}(X) = h(X)$ 
holds for $X\subset V$.  
\finbox 
\medskip

By Lemma \ref{hbase}, Theorem~\ref{THequidotone}
can be applied again to the
general orientation problem covering a non-negative and crossing
supermodular set-function $h$ in the same way as it was applied in the
special case of strong orientation above, but we formulate the result
only for the special case of in-degree constrained and
$k$-edge-connected orientations.

\begin{theorem} 
Let $G=(V,E)$ be an undirected graph endowed with a lower bound function 
$f:V\rightarrow {\bf Z} \cup \{-\infty \}$
and an upper bound function
$g:V\rightarrow {\bf Z} \cup \{+\infty \}$ 
with $f\leq g$.  
A $k$-edge-connected and in-degree constrained orientation $D$ of $G$ 
is decreasingly minimal 
if and only if 
there are no two nodes $s$ and $t$ for which 
$\varrho(t) \geq \varrho(s)+2$,  \
$\varrho(t)>f(t)$, $\varrho(s)<g(s)$, \
and there are $k+1$ arc-disjoint $st$-dipaths.  
\finbox 
\end{theorem}

The theorem can be extended even further to in-degree constrained and
$(k,\ell)$-edge-connected orientations ($\ell\leq k$).

\paragraph{Locally $k$-edge-connected orientations} 
For a subset $S$ of nodes of a digraph $D=(V,A)$, we say that $D$ is {\bf
locally $k$-edge-connected} in $S$ if there are $k$ arc-disjoint
dipaths in $D$ from any node $s$ of $S$ to any other node $t$ of $S$.
This is equivalent to requiring (by Menger's theorem) that $\varrho _D(X)\geq k$
holds for each subset $X\subset V$ for which neither $X\cap S$ nor
$S-X$ is empty.  
Clearly, when $S=V$, we are back at ordinary $k$-edge-connectivity.  
Kir{\'a}ly and Lau \cite{Kiraly-Lau} considered a closely related concept when the
requirement for $D$ is that there are $k$ edge-disjoint dipaths from a
specified root-node $r$ of $S$ to every other node of $S$.  Such a
digraph may be called locally rooted $k$-edge-connected in $S$.

By relying on a base-polyhedral approach, Kir\'aly and Lau
\cite{Kiraly-Lau} proved (Theorem 1.3) that there exists a polynomial
algorithm for finding an orientation of $G$ which is locally rooted
$k$-edge-connected in $S$ {\em and} the orientation meets an in-degree
specification 
$m_{\overline{S}}$ 
on the nodes in $V-S$.  By using an
analogous approach, we can solve the orientation problem in which,
beside the in-degree specification on $V-S$, the requirement is local
$k$-edge-connectivity in $S$.  To this end, define a set-function $h$
on $S$ as follows:
\[
h(Z):= \begin{cases}
 0 & \quad \hbox{if} \quad Z=\emptyset ,
\cr 
\vert E\vert - \widetilde m_{\overline{S}}(V-S) & \quad \hbox{if} \quad Z=S ,
\cr  
\max \{k + i_{G}(Z\cup X) - \widetilde m_{\overline{S}}(X):  X\subseteq V-S \} 
   & \quad \hbox{if} \quad \emptyset \subset Z\subset S. 
\end{cases} 
\]
It is not difficult to see that $h$ is a non-negative, crossing
supermodular function on $S$.  Let $B_{S}:=B'(h)$ denote the integral
base-polyhedron defined by $h$.  
By relying on Lemma \ref{hbase}, one can derive the following.

\begin{theorem}  \label{hmbase} 
Let $G=(V,E)$ be an undirected graph with a specified subset $S$ of $V$.  
Let $m_{\overline{S}}$ be an in-degree specification on $V-S$.  
The set of in-degree vectors of those orientations of $G$
which are locally $k$-edge-connected in $S$ and 
in-degree specified in $V-S$ is an M-convex set $\odotZ{B}$, 
namely, $\odotZ{B} = \{(m_{S},m_{\overline{S}}):  m_{S}\in \odotZ{B}_{S}\}$.  
\finbox 
\end{theorem}

By this theorem, we can determine a decreasingly minimal orientation
among those which are $k$-edge-connected in $S$ and in-degree
specified in $V-S$.  
Even additional in-degree constraints can be imposed on the elements of $S$.

\begin{remark} \rm \label{RMbernath}
One may consider the orientation problem
$(*)$ where, beside local $k$-edge-connectivity in $S$, we impose
upper and lower bounds on the in-degrees of nodes outside $S$ 
(instead of exact values).  
Bern{\'a}th et al.~\cite{Bernath.etal} proved that
the in-degree constrained version of Nash-Willams' well-balanced
orientation problem \cite{NWir} is NP-complete.  
By using a similar method, Bern\'ath \cite{Bernath.oral} recently proved 
that already the problem $(*)$ is NP-complete, even in the special case 
with $k=1$.
\finbox
\end{remark}

\paragraph{Hypergraph orientation}

Let $H=(V,{\cal E})$ be a hypergraph for which we assume that each
hyperedge has at least 2 nodes.  
Orienting a hyperedge $Z$ means that
we designate an element $z$ of $Z$ as its head-node.  
A hyperedge $Z$ with a designated head-node $z\in Z$ 
is a directed hyperedge denoted by $(Z,z)$.  
Orienting a hypergraph means the operation of orienting
each of its hyperedges.  
We say that a directed hyperedge $(Z,z)$
enters a subset $X$ of nodes if $z\in X$ and $Z-X\not =\emptyset $. 
A directed hypergraph is called $k$-edge-connected 
if the in-degree of every non-empty proper subset of nodes is at least $k$.

The following result was proved in \cite{FrankJ50} 
(see, also Theorem 2.22 in the survey paper \cite{FrankS12}).

\begin{theorem}  
The set of in-degree vectors of $k$-edge-connected and
in-degree constrained orientations of a hypergraph forms an M-convex set.  
\finbox 
\end{theorem}  

Therefore we can apply the general results obtained 
for decreasing minimization over M-convex sets.

\paragraph{Acknowledgement} 
We are particularly grateful to an anonymous referee of the paper
whose fundamental suggestions made a significant contribution to
the final form of our work.
This research was supported through the program ``Research in Pairs''
by the Mathematisches Forschungsinstitut Oberwolfach in 2019.
The two weeks we could spend at Oberwolfach provided 
an exceptional opportunity to conduct particularly intensive research.
The research was partially supported by the
National Research, Development and Innovation Fund of Hungary
(FK\_18) -- No. NKFI-128673,
and by CREST, JST, Grant Number JPMJCR14D2, Japan, 
and JSPS KAKENHI Grant Numbers JP26280004, JP20K11697.





\end{document}